# VARIATIONS OF THE SOLUTION TO A STOCHASTIC HEAT EQUATION

By Jason Swanson[1]

*University of Wisconsin–Madison*

We consider the solution to a stochastic heat equation. This solution is a random function of time and space. For a fixed point in space, the resulting random function of time, $F(t)$, has a nontrivial quartic variation. This process, therefore, has infinite quadratic variation and is not a semimartingale. It follows that the classical Itô calculus does not apply. Motivated by heuristic ideas about a possible new calculus for this process, we are led to study modifications of the quadratic variation. Namely, we modify each term in the sum of the squares of the increments so that it has mean zero. We then show that these sums, as functions of $t$, converge weakly to Brownian motion.

**1. Introduction.** Let $u(t,x)$ denote the solution to the stochastic heat equation $u_t = \frac{1}{2} u_{xx} + \dot{W}(t,x)$, with initial conditions $u(0,x) \equiv 0$, where $\dot{W}$ is a space-time white noise on $[0,\infty) \times \mathbb{R}$. That is,

$$u(t,x) = \int_{[0,t] \times \mathbb{R}} p(t-r, x-y) W(dr \times dy),$$

where $p(t,x) = (2\pi t)^{-1/2} e^{-x^2/2t}$ is the heat kernel. Let $F(t) = u(t,x)$, where $x \in \mathbb{R}$ is fixed. In Section 2 we show that $F$ is a centered Gaussian process with covariance function

$$EF(s)F(t) = (2\pi)^{-1/2}(|t+s|^{1/2} - |t-s|^{1/2}),$$

Received December 2005; revised March 2007.

[1]Supported in part by the VIGRE grants of both University of Washington and University of Wisconsin–Madison.

and that $F$ has a nontrivial quartic variation. That is, let $\Pi = \{0 = t_0 < t_1 < t_2 < \cdots\}$, where $t_j \uparrow \infty$, and suppose that $|\Pi| = \sup(t_j - t_{j-1}) < \infty$. If

$$V_\Pi(t) = \sum_{j=1}^{N(t)} |F(t_j) - F(t_{j-1})|^4,$$

where $N(t) = \max\{j : t_j \leq t\}$, then

$$\lim_{|\Pi| \to 0} E\left[\sup_{0 \leq t \leq T} \left|V_\Pi(t) - \frac{6}{\pi} t\right|^2\right] = 0.$$

(See Theorem 2.3.) It follows that $F$ is not a semimartingale, so a stochastic integral with respect to $F$ cannot be defined in the classical Itô sense. (It should be remarked that for a large class of parabolic SPDEs, one obtains better regularity results when the solution $u$ is viewed as a process $t \mapsto u(t, \cdot)$ taking values in Sobolev space, rather than for each fixed $x$. Denis [4] has shown that such processes are in fact Dirichlet processes. Also see Krylov [12].)

In this paper and its sequel, we wish to construct a stochastic integral with respect to $F$ which is a limit of discrete Riemann sums. This construction is based on the heuristic ideas of Chris Burdzy, which were communicated to me during my time as a graduate student. Before elaborating on this construction, it is worth mentioning that, for fixed $x$, the process $t \mapsto u(t, x)$ shares many properties with $B^{1/4}$, the fractional Brownian motion (fBm) with Hurst parameter $H = 1/4$. Several different stochastic integrals with respect to fBm have been developed, and there is a wide literature on this topic. See, for example, Decreusefond [3] and the references therein for a survey of many of these constructions.

We consider discrete Riemann sums over a uniformly spaced time partition $t_j = j\Delta t$, where $\Delta t = n^{-1}$. Let $\Delta F_j = F(t_j) - F(t_{j-1})$. Direct computations with the covariance function demonstrate that

$$E \sum_{j=1}^{\lfloor nt \rfloor} F(t_{j-1}) \Delta F_j \quad \text{and} \quad E \sum_{j=1}^{\lfloor nt \rfloor} F(t_j) \Delta F_j$$

both diverge, showing that left and right endpoint Riemann sums are untenable. We therefore consider Stratonovich-type Riemann sums. There are two kinds of Stratonovich sums that one might consider. The first corresponds to the so-called "trapezoid rule" of elementary calculus and is given by

$$\Phi_T(\Delta t) = \sum_{j=1}^{\lfloor nt \rfloor} \tfrac{1}{2} \left(g'(F(t_{j-1})) + g'(F(t_j))\right) \Delta F_j,$$

where, for now, we take $g$ to be a smooth function. The second corresponds to the so-called "midpoint rule" and is given by

$$\Phi_M(\Delta t) = \sum_{j=1}^{\lfloor nt/2 \rfloor} g'(F(t_{2j-1}))(F(t_{2j}) - F(t_{2j-2})).$$

The midpoint Riemann sum can also be computed using the value of the integrand at points with even index, in which case we would consider

$$\widehat{\Phi}_M(\Delta t) = \sum_{j=1}^{\lfloor nt/2 \rfloor} g'(F(t_{2j}))(F(t_{2j+1}) - F(t_{2j-1})).$$

Note that

$$\Phi_M(\Delta t) + \widehat{\Phi}_M(\Delta t) = \sum_{j=2}^{2\lfloor nt/2 \rfloor+1} g'(F(t_{j-1}))\Delta F_j + \sum_{j=1}^{2\lfloor nt/2 \rfloor} g'(F(t_j))\Delta F_j,$$

so that $\Phi_T(\Delta t) \approx \frac{1}{2}(\Phi_M(\Delta t) + \widehat{\Phi}_M(\Delta t))$, where "$\approx$" means the difference goes to zero uniformly on compacts in probability (ucp) as $\Delta t \to 0$.

One approach to studying these Riemann sums is through a regularization procedure developed by Russo, Vallois and coauthors [7, 8, 14, 15]. For instance, to regularize the trapezoid sum, we define

$$\Phi_T(\Delta t, \varepsilon) = \frac{\Delta t}{2\varepsilon} \sum_{j=0}^{\lfloor nt \rfloor-1} (g'(F(t_j)) + g'(F(t_j + \varepsilon)))(F(t_j + \varepsilon) - F(t_j)),$$

so that $\Phi_T(\Delta t) = \Phi_T(\Delta t, \Delta t)$. We then consider

$$\begin{aligned}
&\lim_{\varepsilon \to 0} \lim_{\Delta t \to 0} \Phi_T(\Delta t, \varepsilon) \\
(1.1) \quad &= \lim_{\varepsilon \to 0} \frac{1}{2\varepsilon} \int_0^t (g'(F(s)) + g'(F(s+\varepsilon)))(F(s+\varepsilon) - F(s))\, ds.
\end{aligned}$$

If this limit exists in probability, it is called the symmetric integral and is denoted by $\int g'(F)\, d^\circ F$. In the case that $F = B^{1/4}$, Gradinaru, Russo and Vallois [8] have shown that the symmetric integral exists, and is simply equal to $g(F(t)) - g(F(0))$. In fact, this result holds for any Hurst parameter $H > 1/6$, which was proven independently by Gradinaru, Nourdin, Russo and Vallois [7] and Cheridito and Nualart [2]. Similarly, we can regularize the midpoint sum by defining

$$\Phi_M(\Delta t, \varepsilon) = \frac{\Delta t}{\varepsilon} \sum_{j=1}^{\lfloor nt/2 \rfloor} g'(F(t_{2j-1}))(F(t_{2j-1} + \varepsilon) - F((t_{2j-1} - \varepsilon) \vee 0)).$$

Again, $\Phi_M(\Delta t) = \Phi_M(\Delta t, \Delta t)$, and this time we find that

(1.2)     $\lim_{\varepsilon \to 0} \lim_{\Delta t \to 0} \Phi_M(\Delta t, \varepsilon)$

$$= \lim_{\varepsilon \to 0} \frac{1}{2\varepsilon} \int_0^t g'(F(s))(F(s+\varepsilon) - F((s-\varepsilon) \vee 0)) \, ds.$$

Using a change of variables, we can see that the right-hand sides of equations (1.1) and (1.2) are equal. In other words, both the trapezoid and the midpoint Riemann sums have the same limit under the regularization procedure.

It is natural to suspect that similar results hold when the regularization procedure is abandoned and we work directly with the discrete Riemann sums. To investigate this, let us consider the Taylor expansion

(1.3)   $g(x+h_1) - g(x+h_2) = \sum_{j=1}^4 \frac{1}{j!} g^{(j)}(x)(h_1^j - h_2^j) + R(x, h_1) - R(x, h_2),$

where

$$R(x, h) = \frac{1}{4!} \int_0^h (h-t)^4 g^{(5)}(x+t) \, dt.$$

In particular, if $M = \sup|g^{(5)}(t)|$, where the supremum is taken over all $t$ between $x$ and $x+h$, then

(1.4)                        $|R(x, h)| \le M|h|^5/5!.$

Substituting $x = F(t_{2j-1})$, $h_1 = \Delta F_{2j}$, and $h_2 = -\Delta F_{2j-1}$ into 1.3, we have

$g(F(t_{2j})) - g(F(t_{2j-2}))$

$\qquad = \sum_{j=1}^4 \frac{1}{j!} g^{(j)}(F(t_{2j-1}))(\Delta F_{2j}^j - (-1)^j \Delta F_{2j-1}^j)$

$\qquad\qquad + R(F(t_{2j-1}), \Delta F_{2j}) - R(F(t_{2j-1}), -\Delta F_{2j-1}).$

Substituting this into the telescoping sum

$$g(F(t)) = g(F(0)) + \sum_{j=1}^N \{g(F(t_{2j})) - g(F(t_{2j-2}))\}$$

$$+ g(F(t)) - g(F(t_{2N})),$$

where $N = \lfloor nt/2 \rfloor$, we have

$\Phi_M(\Delta t) = g(F(t)) - g(F(0))$

(1.5)        $- \frac{1}{2} \sum_{j=1}^N g''(F(t_{2j-1}))(\Delta F_{2j}^2 - \Delta F_{2j-1}^2)$

$$- \tfrac{1}{6} \sum_{j=1}^{N} g'''(F(t_{2j-1}))(\Delta F_{2j}^3 + \Delta F_{2j-1}^3) - \varepsilon_1 - \varepsilon_2 - \varepsilon_3,$$

with

$$\varepsilon_1 = \tfrac{1}{24} \sum_{j=1}^{N} g^{(4)}(F(t_{2j-1}))\Delta F_{2j}^4 - \tfrac{1}{24} \sum_{j=1}^{N} g^{(4)}(F(t_{2j-1}))\Delta F_{2j-1}^4,$$

$$\varepsilon_2 = \sum_{j=1}^{N} R(F(t_{2j-1}), \Delta F_{2j}) - \sum_{j=1}^{N} R(F(t_{2j-1}), -\Delta F_{2j-1}),$$

$$\varepsilon_3 = g(F(t)) - g(F(t_{2N})).$$

By continuity, $\varepsilon_3 \to 0$ ucp. Using Theorem 2.3 and (1.4), we can show that, under suitable assumptions on $g$, $\varepsilon_2 \to 0$ ucp. Similarly, using Theorem 2.3, we can show that both summations in the definition of $\varepsilon_1$ converge to

$$\frac{1}{8\pi} \int_0^t g^{(4)}(F(s)) \, ds,$$

so that $\varepsilon_1 \to 0$ ucp. As a result, we have

$$(1.6) \qquad \begin{aligned} \Phi_M(\Delta t) &\approx g(F(t)) - g(F(0)) - \tfrac{1}{2} \sum_{j=1}^{N} g''(F(t_{2j-1}))(\Delta F_{2j}^2 - \Delta F_{2j-1}^2) \\ &\quad - \tfrac{1}{6} \sum_{j=1}^{N} g'''(F(t_{2j-1}))(\Delta F_{2j}^3 + \Delta F_{2j-1}^3). \end{aligned}$$

By similar reasoning, we also have

$$\begin{aligned} \widehat{\Phi}_M(\Delta t) &\approx g(F(t)) - g(F(0)) - \tfrac{1}{2} \sum_{j=1}^{N} g''(F(t_{2j}))(\Delta F_{2j+1}^2 - \Delta F_{2j}^2) \\ &\quad - \tfrac{1}{6} \sum_{j=1}^{N} g'''(F(t_{2j}))(\Delta F_{2j+1}^3 + \Delta F_{2j}^3). \end{aligned}$$

Taking the average of these two gives

$$\begin{aligned} \Phi_T(\Delta t) &\approx g(F(t)) - g(F(0)) + \tfrac{1}{4} \sum_{j=1}^{\lfloor nt \rfloor} (g''(F(t_j)) - g''(F(t_{j-1})))\Delta F_j^2 \\ &\quad - \tfrac{1}{12} \sum_{j=1}^{\lfloor nt \rfloor} g'''(F(t_j))(\Delta F_{j+1}^3 + \Delta F_j^3). \end{aligned}$$

Using the Taylor expansion $f(b) - f(a) = \frac{1}{2}(f'(a) + f'(b))(b - a) + o(|b-a|^2)$ with $f = g''$, we then have

$$\Phi_T(\Delta t) \approx g(F(t)) - g(F(0)) + \frac{1}{24} \sum_{j=1}^{\lfloor nt \rfloor} g'''(F(t_j))(\Delta F_{j+1}^3 + \Delta F_j^3),$$

which is the discrete analog of the expansion used in Gradinaru et al. [8]. In the sequel to this paper, which will be joint work with Chris Burdzy, we will show that the results of Gradinaru et al. for third-order forward and backward integrals extend to this discrete setting. That is, we will show that

$$\lim_{n \to \infty} \sum_{j=1}^{\lfloor nt \rfloor} g'''(F(t_j))\Delta F_{j+1}^3 = -\lim_{n \to \infty} \sum_{j=1}^{\lfloor nt \rfloor} g'''(F(t_j))\Delta F_j^3$$

$$= -\frac{3}{\pi} \int_0^t g^{(4)}(F(s))\, ds$$

ucp. Hence, in this discrete setting, we obtain the same result as Gradinaru et al. That is, $\Phi_T(\Delta t) \to g(F(t)) - g(F(0))$ ucp.

However, this is only for the trapezoid sum. In the discrete setting, without regularization, the convergence of $\Phi_T(\Delta t)$ no longer implies the convergence of $\Phi_M(\Delta t)$, and the results of Gradinaru et al. do not extend to the discrete midpoint sum. We see from (1.6) that to investigate the convergence of $\Phi_M(\Delta t)$, we must investigate the convergence of

$$(1.7) \qquad \sum_{j=1}^{\lfloor nt/2 \rfloor} g''(F(t_{2j-1}))(\Delta F_{2j}^2 - \Delta F_{2j-1}^2).$$

In Proposition 4.7, we show that, when $g'' \equiv 1$, this sum converges in law to $\kappa B$, where $\kappa$ is an explicit positive constant and $B$ is a standard Brownian motion, independent of $F$. This result suggests that we may define $\int_0^t g'(F(s))\, d^M F(s)$ as the limit, in law, of $\Phi_M(\Delta t)$, and that this integral satisfies the change-of-variables formula

$$(1.8) \qquad g(F(t)) = g(F(0)) + \int_0^t g'(F(s))\, d^M F(s) + \frac{\kappa}{2} \int_0^t g''(F(s))\, dB(s).$$

The emergence of a classical Itô integral as a correction term in this formula shows that, unlike the trapezoid sum, the midpoint sum behaves quite differently in the discrete setting than it does under the regularization procedure. In the case that $g(x) = x^2$, equation (1.8) immediately follows from the results in this paper. (See Corollary 4.8.) The extension of (1.8) to a class of sufficiently smooth functions will be the subject of the sequel.

Note that, when $g'' \equiv 1$, each summand in (1.7) is approximately mean zero and has an approximate variance of $\Delta t$. The convergence of this sum to Brownian motion will therefore follow as a special case of the main result of this paper, Theorem 3.8, which is a Donsker-like invariance principle for processes of the form

$$B_n(t) = \sum_{j=1}^{\lfloor nt \rfloor} \sigma_j^2 h_j(\sigma_j^{-1} \Delta F_j),$$

where $\{h_j\}$ is a sequence of random functions and $\sigma_j^2 = E\Delta F_j^2$. The precise assumptions imposed on the functions $\{h_j\}$ are given in Assumption 3.1. Essentially, the functions $h_j(x)$ must grow no faster than $|x|^2$ and must be chosen so that each of the above summands has mean zero. According to Theorem 3.8, the sequence $\{B_n\}$ converges in law to a Brownian motion, independent of $F$, provided that the variance of the increments of $B_n$ converge. In Section 4 we present several examples where the hypotheses of Theorem 3.8 can be verified by straightforward calculations. Chief among these examples is the case $h_j(x) = (-1)^j(x^2 - 1)$, which gives us Proposition 4.7.

The proof of Theorem 3.8 relies, in part, on the fact that $F$ is a Gaussian process whose increments have covariances which decay polynomially. It should be remarked, however, that this decay rate is too slow for existing mixing results such as those in Herrndorf [9], or existing CLTs such as that in Bulinski [1], to be applicable. Instead, we shall appeal to the additional structure that $F$ possesses. Namely, in the proof of Lemma 3.6, we make significant use of the fact that $F$ has a stochastic integral representation as the convolution of a deterministic kernel against a space-time white noise.

It should be emphasized that the conjectured convergence of $\Phi_M(\Delta t)$ for general functions $g$ is only in law, so that the stochastic calculus which would result from (1.8) would be somewhat different from the usual flavors of stochastic calculus we are used to considering. It is also worth mentioning that the midpoint Riemann sum is not unique in its ability to generate an independent Brownian noise term. Such a term appears, for example, in the study of the asymptotic error for the Euler scheme for SDEs driven by Brownian motion (see Jacod and Protter [10]). An independent Brownian noise term is also generated by the trapezoid Riemann sum when it is applied to fractional Brownian motion $B_{1/6}$. In [7] and [2], it is shown that the symmetric integral $\int (B_{1/6})^2 \, d^\circ B_{1/6}$ does not exist. In fact, the variances of the regularized approximations explode. The same is not true for the discrete trapezoid sums. In fact, for any continuous process $X$,

$$X(t)^3 \approx X(0)^3 + \sum_{j=1}^{\lfloor nt \rfloor} (X(t_j)^3 - X(t_{j-1})^3)$$

$$= X(0)^3 + \sum_{j=1}^{\lfloor nt \rfloor} \Delta X_j (X(t_j)^2 + X(t_j)X(t_{j-1}) + X(t_{j-1})^2)$$

$$= X(0)^3 + \sum_{j=1}^{\lfloor nt \rfloor} \Delta X_j \left( \frac{3}{2}(X(t_j)^2 + X(t_{j-1})^2) - \frac{1}{2}\Delta X_j^2 \right)$$

$$= X(0)^3 + 3 \sum_{j=1}^{\lfloor nt \rfloor} \frac{X(t_j)^2 + X(t_{j-1})^2}{2} \Delta X_j - \frac{1}{2} \sum_{j=1}^{\lfloor nt \rfloor} \Delta X_j^3.$$

Nualart and Ortiz [13] have shown that, for $X = B_{1/6}$, this last sum converges in law to a Brownian motion. This illustrates yet another way in which the discrete approach differs from the regularization method.

**2. The quartic variation of $F$.** Define the Hilbert space $H = L^2(\mathbb{R}^2)$ and construct a centered Gaussian process, $I(h)$, indexed by $h \in H$, such that $E[I(g)I(h)] = \int gh$. Recall that $p(t, x) = (2\pi t)^{-1/2} e^{-x^2/2t}$ and for a fixed pair $(t, x)$, let $h_{tx}(r, y) = 1_{[0,t]}(r)p(t - r, x - y) \in H$. Then

$$(2.1) \qquad F(t) = u(t, x) = \int_{[0,t] \times \mathbb{R}} p(t - r, x - y) W(dr \times dy) = I(h_{tx}).$$

Since $F$ is a centered Gaussian process, its law is determined by its covariance function, which is given in the following lemma. We also derive some needed estimates on the increments of $F$.

LEMMA 2.1.  *For all $s, t \in [0, \infty)$,*

$$(2.2) \qquad EF(s)F(t) = \frac{1}{\sqrt{2\pi}} \left( |t + s|^{1/2} - |t - s|^{1/2} \right).$$

*If $0 \le s < t$, then*

$$(2.3) \qquad \left| E|F(t) - F(s)|^2 - \sqrt{\frac{2(t - s)}{\pi}} \right| \le \frac{1}{t^{3/2}} |t - s|^2.$$

*For fixed $\Delta t > 0$, define $t_j = j\Delta t$, and let $\Delta F_j = F(t_j) - F(t_{j-1})$. If $i, j \in \mathbb{N}$ with $i < j$, then*

$$(2.4) \qquad \left| E[\Delta F_i \Delta F_j] + \sqrt{\frac{\Delta t}{2\pi}} \gamma_{j-i} \right| \le \frac{1}{(t_i + t_j)^{3/2}} \Delta t^2,$$

*where $\gamma_j = 2\sqrt{j} - \sqrt{j - 1} - \sqrt{j + 1}$.*

PROOF.  For (2.2), we may assume $s \leq t$. By (2.1),

$$E[F(s)F(t)] = \int_{\mathbb{R}} \int_0^s p(t-r, x-y)p(s-r, x-y)\, dr\, dy$$

$$= \int_0^s \frac{(2\pi)^{-1}}{\sqrt{(t-r)(s-r)}} \int_{\mathbb{R}} \exp\left\{-\frac{(x-y)^2}{2(t-r)} - \frac{(x-y)^2}{2(s-r)}\right\} dy\, dr$$

$$= \int_0^s \frac{(2\pi)^{-1}}{\sqrt{(t-r)(s-r)}} \int_{\mathbb{R}} \exp\left\{-\frac{(x-y)^2(t+s-2r)}{2(t-r)(s-r)}\right\} dy\, dr.$$

Since $(2\pi)^{-1/2} \int \exp\{-(x-y)^2/2c\}\, dy = \sqrt{c}$, we have

$$E[F(s)F(t)] = \frac{1}{\sqrt{2\pi}} \int_0^s \frac{1}{\sqrt{t+s-2r}}\, dr = \frac{1}{\sqrt{2\pi}}(|t+s|^{1/2} - |t-s|^{1/2}),$$

which verifies the formula.

For (2.3), let $0 \leq s < t$. Then (2.2) implies

$$(2.5) \qquad E|F(t) - F(s)|^2 = \frac{1}{\sqrt{\pi}}(\sqrt{t} + \sqrt{s} - \sqrt{2t+2s} + \sqrt{2t-2s}).$$

Thus,

$$\left| E|F(t) - F(s)|^2 - \sqrt{\frac{2(t-s)}{\pi}} \right| = \frac{1}{\sqrt{\pi}}|\sqrt{t} + \sqrt{s} - \sqrt{2t+2s}|$$

$$= \frac{1}{\sqrt{\pi}}\left| \frac{(\sqrt{t} - \sqrt{s})^2}{\sqrt{t} + \sqrt{s} + \sqrt{2t+2s}} \right|,$$

which gives

$$\left| E|F(t) - F(s)|^2 - \sqrt{\frac{2(t-s)}{\pi}} \right| \leq \frac{(\sqrt{t} - \sqrt{s})^2}{\sqrt{\pi}\,(1+\sqrt{2})\sqrt{t}}$$

$$= \frac{|t-s|^2}{\sqrt{\pi}\,(1+\sqrt{2})\sqrt{t}\,(\sqrt{t} + \sqrt{s})^2}.$$

Hence,

$$(2.6) \qquad \left| E|F(t) - F(s)|^2 - \sqrt{\frac{2(t-s)}{\pi}} \right| \leq \frac{1}{\sqrt{\pi}\,(1+\sqrt{2})t^{3/2}}\,|t-s|^2,$$

which proves (2.3).

Finally, for (2.4), fix $i < j$. Observe that for any $k \geq i$,

$$E[F(t_k)\Delta F_i] = E[F(t_k)F(t_i) - F(t_k)F(t_{i-1})]$$

$$= \frac{1}{\sqrt{2\pi}}(\sqrt{t_k + t_i} - \sqrt{t_k - t_i} - \sqrt{t_k + t_{i-1}} + \sqrt{t_k - t_{i-1}})$$

$$= \sqrt{\frac{\Delta t}{2\pi}}(\sqrt{k+i} - \sqrt{k-i} - \sqrt{k+i-1} + \sqrt{k-i+1}).$$

Thus,

$$E[\Delta F_i \Delta F_j] = E[F(t_j)\Delta F_i] - E[F(t_{j-1})\Delta F_i]$$
$$= \sqrt{\frac{\Delta t}{2\pi}}(\sqrt{j+i} - \sqrt{j-i} - \sqrt{j+i-1} + \sqrt{j-i+1}$$
$$- \sqrt{j+i-1} + \sqrt{j-i-1} + \sqrt{j+i-2} - \sqrt{j-i}),$$

which simplifies to

$$(2.7) \qquad E[\Delta F_i \Delta F_j] = -\sqrt{\frac{\Delta t}{2\pi}}(\gamma_{j+i-1} + \gamma_{j-i}).$$

The strict concavity of $x \mapsto \sqrt{x}$ implies that $\gamma_k > 0$ for all $k \in \mathbb{N}$. Also, if we write $\gamma_k = f(k-1) - f(k)$, where $f(x) = \sqrt{x+1} - \sqrt{x}$, then for each $k \geq 2$, the mean value theorem gives $\gamma_k = |f'(k-\theta)|$ for some $\theta \in [0,1]$. Since $|f'(x)| \leq x^{-3/2}/4$, we can easily verify that for all $k \in \mathbb{N}$,

$$(2.8) \qquad 0 < \gamma_k \leq \frac{1}{\sqrt{2}\,k^{3/2}}.$$

Since $j+i-1 \geq (j+i)/2$, we have

$$\left| E[\Delta F_i \Delta F_j] + \sqrt{\frac{\Delta t}{2\pi}}\,\gamma_{j-i} \right| \leq \sqrt{\frac{\Delta t}{2\pi}}\,\frac{1}{\sqrt{2}\,((j+i)/2)^{3/2}} = \sqrt{\frac{2\Delta t}{\pi}}\,\frac{1}{(i+j)^{3/2}},$$

and this proves (2.4).  □

By (2.2), the law of $F(t) = u(t,x)$ does not depend on $x$. We will therefore assume that $x = 0$. Note that (2.6) implies

$$(2.9) \qquad \pi^{-1/2}\sqrt{\Delta t} \leq E\Delta F_j^2 \leq 2\sqrt{\Delta t}$$

for all $j \geq 1$. In particular, since $F$ is Gaussian, we have $E|F(t) - F(s)|^{4n} \leq C_n|t-s|^n$ for all $n$. By the Kolmogorov–Čentsov theorem (see, e.g. Theorem 2.2.8 in [11]), $F$ has a modification which is locally Hölder continuous with exponent $\gamma$ for all $\gamma \in (0, 1/4)$. We will henceforth assume that we are working with such a modification.  □

Also note that (2.7) and (2.8) together imply

$$-\frac{2\Delta t^2}{(t_j - t_i)^{3/2}} = -\frac{2\sqrt{\Delta t}}{(j-i)^{3/2}} \leq E[\Delta F_i \Delta F_j] < 0$$

for all $1 \leq i < j$. In other words, the increments of $F$ are negatively correlated and we have a polynomial bound on the rate of decay of this correlation.

For future reference, let us combine these results into the following single inequality: for all $i, j \in \mathbb{N}$,

$$(2.10) \qquad |E[\Delta F_i \Delta F_j]| \leq \frac{2\sqrt{\Delta t}}{|i - j|^{\sim 3/2}},$$

where we have adopted the notation $x^{\sim r} = (x \vee 1)^r$. In fact, with a little more work, we have the following general result.

LEMMA 2.2. *For all* $0 \leq s < t \leq u < v,$

$$(2.11) \qquad |E[(F(v) - F(u))(F(t) - F(s))]| \leq \sqrt{\frac{2}{\pi}} \frac{|t - s||v - u|}{|u - s|\sqrt{v - t}}.$$

PROOF. Fix $0 \leq s < t$. For any $r > t$, define

$$f(r) = E[F(r)(F(t) - F(s))] = \frac{1}{\sqrt{2\pi}} (\sqrt{r + t} - \sqrt{r - t} - \sqrt{r + s} + \sqrt{r - s}).$$

Then

$$f'(r) = \frac{1}{2\sqrt{2\pi}} \left( \frac{\sqrt{r + s} - \sqrt{r + t}}{\sqrt{(r + t)(r + s)}} - \frac{\sqrt{r - s} - \sqrt{r - t}}{\sqrt{(r - t)(r - s)}} \right).$$

Since

$$\left| \frac{\sqrt{r \pm s} - \sqrt{r \pm t}}{\sqrt{(r \pm t)(r \pm s)}} \right| \leq \frac{|t - s|}{\sqrt{r - t}|r - s|},$$

we have

$$|E[(F(v) - F(u))(F(t) - F(s))]| \leq \frac{1}{\sqrt{2\pi}} \int_u^v \frac{|t - s|}{\sqrt{r - t}|r - s|} \, dr$$
$$\leq \frac{1}{\sqrt{2\pi}} \frac{|t - s|}{|u - s|} \int_u^v \frac{1}{\sqrt{r - t}} \, dr$$
$$= \sqrt{\frac{2}{\pi}} \frac{|t - s|}{|u - s|} (\sqrt{v - t} - \sqrt{u - t})$$
$$\leq \sqrt{\frac{2}{\pi}} \frac{|t - s||v - u|}{|u - s|\sqrt{v - t}}$$

whenever $0 \leq s < t \leq u < v$. $\quad \square$

THEOREM 2.3. *Let*

$$V_\Pi(t) = \sum_{j=1}^{N(t)} |F(t_j) - F(t_{j-1})|^4,$$

where $\Pi = \{0 = t_0 < t_1 < t_2 < \cdots\}$ is a partition of $[0, \infty)$, that is, $t_j \uparrow \infty$ and $N(t) = \max\{j : t_j \leq t\}$. Let $|\Pi| = \sup(t_j - t_{j-1}) < \infty$. Then

$$\lim_{|\Pi| \to 0} E \left[ \sup_{0 \leq t \leq T} \left| V_\Pi(t) - \frac{6}{\pi} t \right|^2 \right] = 0$$

for all $T > 0$.

PROOF. Since $V_\Pi$ is monotone, it will suffice to show that $V_\Pi(t) \to 6t/\pi$ in $L^2$ for each fixed $t$. In what follows, $C$ is a finite, positive constant that may change value from line to line. Fix $t \geq 0$ and let $N = N(t)$. For each $j$, let $\Delta F_j = F(t_j) - F(t_{j-1})$, $\sigma_j^2 = E\Delta F_j^2$, and $\Delta t_j = t_j - t_{j-1}$. Note that

$$V_\Pi(t) = \sum_{j=1}^N (\Delta F_j^4 - 3\sigma_j^4) + \sum_{j=1}^N \left( 3\sigma_j^4 - \frac{6}{\pi} \Delta t_j \right) + \frac{6}{\pi}(t_N - t).$$

By (2.3),

$$\left| 3\sigma_j^4 - \frac{6}{\pi} \Delta t_j \right| = 3 \left| \sigma_j^2 + \sqrt{\frac{2\Delta t_j}{\pi}} \right| \left| \sigma_j^2 - \sqrt{\frac{2\Delta t_j}{\pi}} \right| \leq \frac{C}{t_j^{3/2}} \Delta t_j^{5/2} \leq \frac{C}{t_j^{3/4}} \Delta t_j^{7/4}.$$

Thus,

$$(2.12) \qquad \left| \sum_{j=1}^N \left( 3\sigma_j^4 - \frac{6}{\pi} \Delta t_j \right) \right| \leq C |\Pi|^{3/4} \sum_{j=1}^N \frac{\Delta t_j}{t_j^{3/4}},$$

which tends to zero as $|\Pi| \to 0$ since $\int_0^t x^{-3/4} \, dx < \infty$.

To complete the proof, we will need the following fact about Gaussian random variables. Let $X_1, X_2$ be mean zero, jointly normal random variables with variances $\sigma_j^2$. If $\rho = (\sigma_1 \sigma_2)^{-1} E[X_1 X_2]$, then

$$(2.13) \qquad E[X_1^4 X_2^4] = \sigma_1^4 \sigma_2^4 (24\rho^4 + 72\rho^2 + 9).$$

Applying this in our context, let $\rho_{ij} = (\sigma_i \sigma_j)^{-1} E[\Delta F_i \Delta F_j]$ and write

$$E \left| \sum_{j=1}^N (\Delta F_j^4 - 3\sigma_j^4) \right|^2 \leq \sum_{i=1}^N \sum_{j=1}^N |E[(\Delta F_i^4 - 3\sigma_i^4)(\Delta F_j^4 - 3\sigma_j^4)]|$$

$$= \sum_{i=1}^N \sum_{j=1}^N |E[\Delta F_i^4 \Delta F_j^4] - 9\sigma_i^4 \sigma_j^4|.$$

Then by (2.13), we have

$$E \left| \sum_{j=1}^N (\Delta F_j^4 - 3\sigma_j^4) \right|^2 \leq C \sum_{i=1}^N \sum_{j=1}^N \sigma_i^4 \sigma_j^4 \rho_{ij}^2$$

$$= C \sum_{i=1}^{N} \sum_{j=1}^{N} \sigma_i^2 \sigma_j^2 |E[\Delta F_i \Delta F_j]|^2$$

$$\leq C \sum_{i=1}^{N} \sum_{j=1}^{N} \Delta t_i^{1/2} \Delta t_j^{1/2} |E[\Delta F_i \Delta F_j]|^2.$$

By Hölder's inequality, $|E[\Delta F_i \Delta F_j]|^2 \leq \Delta t_i^{1/2} \Delta t_j^{1/2}$, so it will suffice to show that

$$\sum_{i=1}^{N-2} \sum_{j=i+2}^{N} \Delta t_i^{1/2} \Delta t_j^{1/2} |E[\Delta F_i \Delta F_j]|^2 \to 0$$

as $|\Pi| \to 0$. For this, suppose $j > i + 1$. By (2.11),

$$|E[\Delta F_i \Delta F_j]|^2 \leq \frac{C \Delta t_i^2 \Delta t_j^2}{|t_{j-1} - t_{i-1}|^2 |t_j - t_i|}$$

$$\leq \frac{C \Delta t_i^{1/2} \Delta t_j^{5/4}}{|t_{j-1} - t_{i-1}|^{1/2} |t_j - t_i|^{1/4}}$$

$$\leq C |\Pi|^{3/4} \frac{\Delta t_i^{1/2} \Delta t_j^{1/2}}{|t_{j-1} - t_i|^{3/4}}.$$

Hence,

$$\sum_{i=1}^{N-2} \sum_{j=i+2}^{N} \Delta t_i^{1/2} \Delta t_j^{1/2} |E[\Delta F_i \Delta F_j]|^2$$

$$\leq C |\Pi|^{3/4} \sum_{i=1}^{N-2} \sum_{j=i+2}^{N} |t_{j-1} - t_i|^{-3/4} \Delta t_i \Delta t_j,$$

which tends to zero as $|\Pi| \to 0$ since $\int_0^t \int_0^t |x - y|^{-3/4} \, dx \, dy < \infty$.   $\square$

**3. Main result.** Let us now specialize to the uniform partition. That is, for fixed $n \in \mathbb{N}$, let $\Delta t = n^{-1}$, $t_j = j \Delta t$ and $\Delta F_j = F(t_j) - F(t_{j-1})$. We wish to consider sums of the form $\sum_{j=1}^{\lfloor nt \rfloor} g_j(\Delta F_j)$, where $\{g_j\}$ is a sequence of random functions. We will write these functions in the form $g_j(x) = \sigma_j^2 h_j(\sigma_j^{-1} x)$, where $\sigma_j^2 = E \Delta F_j^2$.

ASSUMPTION 3.1. Let $\{h_j(x) : x \in \mathbb{R}\}$ be a sequence of independent stochastic processes which are almost surely continuously differentiable. Assume there exists a constant $L$ such that $E h_j(0)^2 \leq L$ and $E h_j'(0)^2 \leq L$ for all $j$. Also assume that for each $j$,

$$|h_j'(x) - h_j'(y)| \leq L_j |x - y| \tag{3.1}$$

for all $x, y \in \mathbb{R}$, where $EL_j^2 \le L$. Finally, assume that

$$(3.2) \qquad Eh_j(X) = 0$$

and

$$(3.3) \qquad |Eh_i(X)h_j(Y)| \le L|\rho|$$

whenever $X$ and $Y$ which are independent of $\{h_j\}$ and are jointly normal with mean zero, variance one, and covariance $\rho = EXY$.

REMARK 3.2.   We may assume that each $L_j$ is $\sigma(h_j)$-measurable. In particular, $\{L_j\}$ is a sequence of independent random variables. Also, since $Eh'_j(0)^2 \le L$, we may assume that

$$(3.4) \qquad |h'_j(x)| \le L_j(1 + |x|)$$

for all $j$. Similarly, since $Eh_j(0)^2 \le L$, we may assume that

$$(3.5) \qquad |h_j(x)| \le L_j(1 + |x|^2)$$

for all $j$.

LEMMA 3.3.   *Let* $\{h_j\}$ *satisfy Assumption 3.1. Let* $X_1, \ldots, X_4$ *be mean zero, jointly normal random variables, independent of the sequence* $\{h_j\}$, *such that* $EX_j^2 = 1$ *and* $\rho_{ij} = EX_iX_j$. *Then there exists a finite constant* $C$, *that depends only on* $L$, *such that*

$$(3.6) \qquad \left| E\prod_{j=1}^{4} h_j(X_j) \right| \le C\left( |\rho_{12}\rho_{34}| + \frac{1}{\sqrt{1 - \rho_{12}^2}} \max_{i \le 2 < j} |\rho_{ij}| \right)$$

*whenever* $|\rho_{12}| < 1$. *Moreover,*

$$(3.7) \qquad \left| E\prod_{j=1}^{4} h_j(X_j) \right| \le C \max_{2 \le j \le 4} |\rho_{1j}|.$$

*Furthermore, there exists* $\varepsilon > 0$ *such that*

$$(3.8) \qquad \left| E\prod_{j=1}^{4} h_j(X_j) \right| \le CM^2$$

*whenever* $M = \max\{|\rho_{ij}| : i \ne j\} < \varepsilon$.

PROOF.   In the proofs in this section, $C$ will denote a finite, positive constant that depends only on $L$, which may change value from line to line. Let us first record some observations. By (3.4), with probability one,

$$(3.9) \qquad \begin{aligned} |h_j(y) - h_j(x)| &= \left| \int_0^1 (y - x)h'_j(x + t(y - x))\, dt \right| \\ &\le L_j|y - x|(1 + |x| + |y - x|). \end{aligned}$$

Also,

$$h_j(y) - h_j(x) - (y - x)h_j'(x)$$
$$= \int_0^1 (y - x)(h_j'(x + t(y - x)) - h_j'(x))\,dt,$$

so by (3.1),

$$(3.10) \qquad |h_j(y) - h_j(x) - (y - x)h_j'(x)| \le |y - x| \int_0^1 L_j t |y - x|\,dt$$
$$\le L_j |y - x|^2$$

for all $x, y \in \mathbb{R}$. Also, by (3.4) and (3.5), if we define a stochastic process on $\mathbb{R}^n$ by $G(x) = \prod_{j=1}^n h_j(x_j)$, then $G$ is almost surely continuously differentiable with $|\partial_j G(x)| \le C(\prod_{j=1}^n L_j)(1 + |x|^{2n-1})$. Hence,

$$(3.11) \qquad |G(y) - G(x)| = \left| \int_0^t \frac{d}{dt} G(x + t(y - x))\,dt \right|$$
$$= \left| \int_0^t (y - x) \cdot \nabla G(x + t(y - x))\,dt \right|$$
$$\le C\left( \prod_{j=1}^n L_j \right) n |y - x| (1 + |x|^{2n-1} + |y - x|^{2n-1})$$

for all $x, y \in \mathbb{R}^n$.

Now let $Y_1 = (X_1, X_2)^T$ and $Y_2 = (X_3, X_4)^T$. If $|\rho_{12}| < 1$, then we may define the matrix $A = (EY_2 Y_1^T)(EY_1 Y_1^T)^{-1}$. Note that

$$(3.12) \qquad |A| \le \frac{C}{\sqrt{1 - \rho_{12}^2}} \max_{i \le 2 < j} |\rho_{ij}|.$$

Let $\bar{Y}_2 = Y_2 - AY_1$, so that $E\bar{Y}_2 Y_1^T = 0$, which implies $\bar{Y}_2$ and $Y_1$ are independent, and define stochastic processes on $\mathbb{R}^2$ by

$$F_{ij}(x) = h_i(x_1)h_j(x_2),$$

so that $\prod_{j=1}^4 h_j(X_j) = F_{12}(Y_1) F_{34}(Y_2)$.

Also define $\bar{X} = (X_2, X_3, X_4)^T$ and $c = (\rho_{12}, \rho_{13}, \rho_{14})^T$. Note that $X_1$ and $\bar{X} - cX_1$ are independent. Define a process on $\mathbb{R}^3$ by

$$F(x) = h_2(x_1)h_3(x_2)h_4(x_3),$$

so that $\prod_{j=1}^4 h_j(X_j) = h_1(X_1)F(\bar{X})$.

Let $\Sigma = E\bar{X}\bar{X}^T$. If $M$ is sufficiently small, then $\Sigma$ is invertible, we may define $a = \Sigma^{-1}c$, and we have $|a| \le CM$ and $|a^T c| \le CM^2 < 3/4$. Note that $a^T c = E|a^T \bar{X}|^2 \ge 0$. Let $\sigma = (1 - a^T c)^{-1/2}$ so that $1 \le \sigma < 2$ and

$$(3.13) \qquad \sigma - 1 = \frac{1 - \sigma^{-1}}{\sigma^{-1}} < \frac{2a^T c}{1 + \sigma^{-1}} < 2a^T c.$$

Define $U = \sigma(X_1 - a^T\bar{X})$. Note that

$$E[U\bar{X}^T] = \sigma(c^T - a^T\Sigma) = 0,$$

so that $U$ and $\bar{X}$ are independent. Hence,

$$\sigma^2 = E(\sigma X_1)^2 = EU^2 + \sigma^2 E|a^T\bar{X}|^2 = EU^2 + \sigma^2 a^T c,$$

so that $U$ is normal with mean zero and variance one.

For the proof of (3.6), we have

$$E\prod_{j=1}^4 h_j(X_j) = E[F_{12}(Y_1)F_{34}(\bar{Y}_2)] + E[F_{12}(Y_1)(F_{34}(Y_2) - F_{34}(\bar{Y}_2))]$$

$$(3.14)$$
$$= EF_{12}(Y_1)EF_{34}(\bar{Y}_2) + E[F_{12}(Y_1)(F_{34}(Y_2) - F_{34}(\bar{Y}_2))]$$

and

$$(3.15) \qquad EF_{34}(\bar{Y}_2) = EF_{34}(Y_2) - E[F_{34}(Y_2) - F_{34}(\bar{Y}_2)].$$

By (3.11),

$$|F_{34}(Y_2) - F_{34}(\bar{Y}_2)| \le CL_3 L_4 |AY_1|(1 + |Y_2|^3 + |AY_1|^3).$$

Note that $EL_3^2 L_4^2 = EL_3^2 EL_4^2 \le L^2$. Also, since $E|Y_2|^2 = E|\bar{Y}_2|^2 + E|AY_1|^2$, we see that the components of $AY_1$ are jointly normal with mean zero and a variance which is bounded by a constant independent of $\{\rho_{ij}\}$. Hence, Hölder's inequality gives

$$E|F_{34}(Y_2) - F_{34}(\bar{Y}_2)|^2 \le C|A|^2(E|Y_1|^4)^{1/2}(1 + E|Y_2|^{12} + E|AY_1|^{12})^{1/2}$$
$$\le C|A|^2.$$

By (3.12),

$$(3.16) \qquad (E|F_{34}(Y_2) - F_{34}(\bar{Y}_2)|^2)^{1/2} \le \frac{C}{\sqrt{1 - \rho_{12}^2}} \max_{i \le 2 < j} |\rho_{ij}|.$$

Hence, by (3.15),

$$(3.17) \qquad |EF_{34}(\bar{Y}_2)| \le |EF_{34}(Y_2)| + \frac{C}{\sqrt{1 - \rho_{12}^2}} \max_{i \le 2 < j} |\rho_{ij}|.$$

Note that (3.5) implies $E|F_{12}(Y_1)|^2 \le C$. Therefore, using (3.14), (3.15), (3.16) and Hölder's inequality, we have

$$\left| E\prod_{j=1}^4 h_j(X_j) \right| \le |EF_{12}(Y_1)EF_{34}(Y_2)| + \frac{C}{\sqrt{1 - \rho_{12}^2}} \max_{i \le 2 < j} |\rho_{ij}|.$$

By (3.3), this completes the proof of (3.6).

For (3.7), we have

$$E \prod_{j=1}^{4} h_j(X_j) = E h_1(X_1) E F(\bar{X} - cX_1)$$
$$+ E[h_1(X_1)(F(\bar{X}) - F(\bar{X} - cX_1))].$$

Since $X_1$ and $\bar{X} - cX_1$ are independent, (3.2) gives

(3.18)
$$E \prod_{j=1}^{4} h_j(X_j) = E[h_1(X_1)(F(\bar{X}) - F(\bar{X} - cX_1))].$$

By Hölder's inequality and (3.5),

$$\left| E \prod_{j=1}^{4} h_j(X_j) \right| \leq C(E|F(\bar{X}) - F(\bar{X} - cX_1)|^2)^{1/2}.$$

By (3.11),

$$|F(\bar{X}) - F(\bar{X} - cX_1)| \leq C\left(\prod_{j=2}^{4} L_j\right) |cX_1|(1 + |\bar{X}|^5 + |cX_1|^5).$$

Hence,

$$E|F(\bar{X}) - F(\bar{X} - cX_1)|^2 \leq C|c|^2,$$

which gives

$$\left| E \prod_{j=1}^{4} h_j(X_j) \right| \leq C|c|,$$

and proves (3.7).

Finally, for (3.8), we begin with an auxiliary result. Note that

$$E\left[ X_2 \prod_{j=2}^{4} h_j(X_j) \right] = E[X_2 h_2(X_2)] E F_{34}(\bar{Y}_2)$$
$$+ E[X_2 h_2(X_2)(F_{34}(Y_2) - F_{34}(\bar{Y}_2))].$$

By Hölder's inequality and (3.5),

$$\left| E\left[ X_2 \prod_{j=2}^{4} h_j(X_j) \right] \right| \leq C|E F_{34}(\bar{Y}_2)| + C(E|F_{34}(Y_2) - F_{34}(\bar{Y}_2)|^2)^{1/2}.$$

If $M$ is sufficiently small, then $|\rho_{12}| \leq C < 1$. Hence, by (3.16), (3.17), and (3.3),

$$\left| E\left[ X_2 \prod_{j=2}^{4} h_j(X_j) \right] \right| \leq CM.$$

It now follows by symmetry that

$$(3.19) \qquad \left| E\left[ v^T \bar{X} \prod_{j=2}^{4} h_j(X_j) \right] \right| \le C|v|M$$

for any $v \in \mathbb{R}^3$.

Returning to the proof of (3.8), since (3.2) implies $Eh_1(U) = 0$ and $U$ and $\bar{X}$ are independent, we have

$$E \prod_{j=1}^{4} h_j(X_j) = E[(h_1(X_1) - h_1(U) - (X_1 - U)h_1'(U))F(\bar{X})]$$

$$+ E[(X_1 - U)h_1'(U)F(\bar{X})].$$

By (3.10),

$$|h_1(X_1) - h_1(U) - (X_1 - U)h_1'(U)| \le L_1|X_1 - U|^2.$$

By (3.13), $|1 - \sigma| \le Ca^Tc \le C|a|M$, so that

$$|X_1 - U| = |(1 - \sigma)X_1 + a^T\bar{X}| \le C|a|(M|X_1| + |\bar{X}|).$$

Hence, using Hölder's inequality and (3.5), we have

$$(3.20) \quad \begin{aligned} \left| E \prod_{j=1}^{4} h_j(X_j) \right| &\le C(E|X_1 - U|^4)^{1/2} + |E[(X_1 - U)h_1'(U)F(\bar{X})]| \\ &\le C|a|^2 + |E[(X_1 - U)h_1'(U)F(\bar{X})]|. \end{aligned}$$

To estimate the second term, note that

$$E[(X_1 - U)h_1'(U)F(\bar{X})] = (1 - \sigma)E[X_1 h_1'(U)F(\bar{X})] + Eh_1'(U) \cdot E[a^T\bar{X}F(\bar{X})].$$

Therefore, by (3.4), (3.5), (3.13) and (3.19),

$$|E[(X_1 - U)h_1'(U)F(\bar{X})]| \le C|1 - \sigma| + C|a|M \le C|a|M.$$

Combining this with (3.20) and recalling that $|a| \le CM$ completes the proof of (3.8). $\quad\square$

COROLLARY 3.4. *Let $\{h_j\}$ be independent of $F$ and satisfy Assumption 3.1. For $k \in \mathbb{N}^4$ with $k_1 \le \cdots \le k_4$, define*

$$(3.21) \qquad \Delta_k = \prod_{j=1}^{4} \sigma_{k_j}^2 h_{k_j}(\sigma_{k_j}^{-1}\Delta F_{k_j}),$$

*where $\sigma_j^2 = E\Delta F_j^2$. Let $x^{\sim r} = (x \vee 1)^r$. Then there exists a finite constant $C$ such that*

$$(3.22) \qquad |E\Delta_k| \le \frac{C\Delta t^2}{(k_4 - k_3)^{\sim 3/2}} \quad and \quad |E\Delta_k| \le \frac{C\Delta t^2}{(k_2 - k_1)^{\sim 3/2}}.$$

*Moreover,*

$$(3.23) \qquad |E\Delta_k| \le C\Big(\frac{1}{(k_4 - k_3)^{\sim 3/2}(k_2 - k_1)^{\sim 3/2}} + \frac{1}{(k_3 - k_2)^{\sim 3/2}}\Big)\Delta t^2$$

*and*

$$(3.24) \qquad\qquad |E\Delta_k| \le \frac{C\Delta t^2}{m^{\sim 3}},$$

*where $m = \min\{k_{i+1} - k_i : 1 \le i < 4\}$.*

PROOF.  Let $X_j = \sigma_{k_j}^{-1}\Delta F_{k_j}$. By (2.9) and (2.10), we have

$$|\rho_{ij}| = |E[X_i X_j]| = \sigma_{k_i}^{-1}\sigma_{k_j}^{-1}|E[\Delta F_{k_i}\Delta F_{k_j}]| \le \frac{2\sqrt{\pi}}{|k_i - k_j|^{\sim 3/2}}.$$

Also,

$$|E\Delta_k| = \Big(\prod_{j=1}^{4}\sigma_{k_j}^2\Big)\Big|E\Big[\prod_{j=1}^{4}h_{k_j}(X_j)\Big]\Big|.$$

This, together with (3.7) and symmetry, yields (3.22).

For (3.23), first note that Hölder's inequality and (3.5) give the trivial bound $|E\Delta_k| \le C\Delta t^2$. Hence, we may assume that at least one of $k_4 - k_3$ and $k_2 - k_1$ is large. Specifically, by symmetry, we may assume that $k_2 - k_1 \ge 4$. In this case, $|\rho_{12}| \le \sqrt{\pi}/4 < 1$. Hence, (3.6) gives

$$\Big(\prod_{j=1}^{4}\sigma_{k_j}^2\Big)\Big|E\Big[\prod_{j=1}^{4}h_{k_j}(X_j)\Big]\Big| \le C\Big(\prod_{j=1}^{4}\sigma_{k_j}^2\Big)\Big(|\rho_{12}\rho_{34}| + \frac{1}{\sqrt{1 - \rho_{12}^2}}\max_{i \le 2 < j}|\rho_{ij}|\Big)$$

and (3.23) is immediate.

As above, we may assume in proving (3.24) that $m$ is large. Therefore, we can assume that $M = \max\{|\rho_{ij}| : i \ne j\} < \varepsilon$. Hence, (3.8) implies

$$\Big(\prod_{j=1}^{4}\sigma_{k_j}^2\Big)\Big|E\Big[\prod_{j=1}^{4}h_{k_j}(X_j)\Big]\Big| \le C\Big(\prod_{j=1}^{4}\sigma_{k_j}^2\Big)M^2,$$

which proves (3.24).  □

PROPOSITION 3.5.  *With notation as in Corollary 3.4, let*

$$(3.25) \qquad\qquad B_n(t) = \sum_{j=1}^{\lfloor nt \rfloor}\sigma_j^2 h_j(\sigma_j^{-1}\Delta F_j).$$

If $\{h_j\}$ is independent of $F$ and satisfies Assumption *3.1,* then there exists a constant $C$ such that

$$E|B_n(t) - B_n(s)|^4 \le C\left(\frac{\lfloor nt\rfloor - \lfloor ns\rfloor}{n}\right)^2 \tag{3.26}$$

for all $0 \le s < t$ and all $n \in \mathbb{N}$. The sequence $\{B_n\}$ is therefore relatively compact in the Skorohod space $D_{\mathbb{R}}[0, \infty)$.

PROOF. To prove (3.26), observe that

$$E|B_n(t) - B_n(s)|^4 = E\left|\sum_{j=ns+1}^{\lfloor nt\rfloor} \sigma_j^2 h_j(\sigma_j^{-1}\Delta F_j)\right|^4.$$

Let

$$S = \{k \in \mathbb{N}^4 \colon \lfloor ns\rfloor + 1 \le k_1 \le \cdots \le k_4 \le \lfloor nt\rfloor\}.$$

For $k \in S$, define $h_i = k_{i+1} - k_i$ and let

$$M = M(k) = \max(h_1, h_2, h_3),$$
$$m = m(k) = \min(h_1, h_2, h_3),$$
$$c = c(k) = \text{med}(h_1, h_2, h_3),$$

where "med" denotes the median function. For $i \in \{1, 2, 3\}$, let $S_i = \{k \in S \colon h_i = M\}$. Define $N = \lfloor nt\rfloor - (\lfloor ns\rfloor + 1)$ and for $j \in \{0, 1, \ldots, N\}$, let $S_i^j = \{k \in S_i \colon M = j\}$. Further define $T_i^\ell = T_i^{j,\ell} = \{k \in S_i^j \colon m = \ell\}$ and $V_i^\nu = V_i^{j,\ell,\nu} = \{k \in T_i^\ell \colon c = \nu\}$.

Recalling (3.21), we now have

$$E\left|\sum_{j=\lfloor ns\rfloor+1}^{\lfloor nt\rfloor} \sigma_j^2 h_j(\sigma_j^{-1}\Delta F_j)\right|^4 \le 4! \sum_{k \in S} |E\Delta_k| \le 4! \sum_{i=1}^{3} \sum_{k \in S_i} |E\Delta_k|. \tag{3.27}$$

Observe that

$$\sum_{k \in S_i} |E\Delta_k| = \sum_{j=0}^{N} \sum_{k \in S_i^j} |E\Delta_k| \tag{3.28}$$

and

$$\sum_{k \in S_i^j} |E\Delta_k| = \sum_{\ell=0}^{\lfloor\sqrt{j}\rfloor} \sum_{k \in T_i^\ell} |E\Delta_k| + \sum_{\ell=\lfloor\sqrt{j}\rfloor+1}^{j} \sum_{k \in T_i^\ell} |E\Delta_k|. \tag{3.29}$$

Begin by considering the first summation. Suppose $0 \le \ell \le \lfloor\sqrt{j}\rfloor$ and write

$$\sum_{k \in T_i^\ell} |E\Delta_k| = \sum_{\nu=\ell}^{j} \sum_{k \in V_i^\nu} |E\Delta_k|.$$

Fix $\nu$ and let $k \in V_i^\nu$ be arbitrary. If $i = 1$, then $j = M = h_1 = k_2 - k_1$. If $i = 3$, then $j = M = h_3 = k_4 - k_3$. In either case, (3.22) gives

$$|E\Delta_k| \le C \frac{1}{j^{\sim 3/2}} \Delta t^2 \le C \left( \frac{1}{(\ell \nu)^{\sim 3/2}} + \frac{1}{j^{\sim 3/2}} \right) \Delta t^2.$$

If $i = 2$, then $j = M = h_2 = k_3 - k_2$ and $\ell \nu = h_3 h_1 = (k_4 - k_3)(k_2 - k_1)$. Hence, by (3.23),

$$|E\Delta_k| \le C \left( \frac{1}{(\ell \nu)^{\sim 3/2}} + \frac{1}{j^{\sim 3/2}} \right) \Delta t^2.$$

Now choose $i' \ne i$ such that $h_{i'} = \ell$. With $i'$ given, $k$ is determined by $k_i$. Since there are two possibilities for $i'$ and $N+1$ possibilities for $k_i$, $|V_i^\nu| \le 2(N+1)$. Therefore,

$$\sum_{\ell=0}^{\lfloor \sqrt{j} \rfloor} \sum_{k \in T_i^\ell} |E\Delta_k| \le C(N+1) \sum_{\ell=0}^{\lfloor \sqrt{j} \rfloor} \sum_{\nu=\ell}^{j} \left( \frac{1}{(\ell \nu)^{\sim 3/2}} + \frac{1}{j^{\sim 3/2}} \right) \Delta t^2$$

$$\le C(N+1) \sum_{\ell=0}^{\lfloor \sqrt{j} \rfloor} \left( \frac{1}{\ell^{\sim 3/2}} + \frac{1}{j^{\sim 1/2}} \right) \Delta t^2$$

$$\le C(N+1) \Delta t^2.$$

For the second summation, suppose $\lfloor \sqrt{j} \rfloor + 1 \le \ell \le j$. (In particular, $j \ge 1$.) In this case, if $k \in T_i^\ell$, then $\ell = m = \min\{k_{i+1} - k_i : 1 \le i < 4\}$, so that by (3.24),

$$|E\Delta_k| \le C \frac{1}{\ell^{\sim 3}} \Delta t^2.$$

Since $|T_i^\ell| = \sum_{\nu=\ell}^{j} |V_i^\nu| \le 2(N+1)j$, we have

$$\sum_{\ell=\sqrt{j}+1}^{j} \sum_{k \in T_i^\ell} |E\Delta_k| \le C(N+1)j \sum_{\ell=\lfloor \sqrt{j} \rfloor + 1}^{j} \frac{1}{\ell^{\sim 3}} \Delta t^2$$

$$\le C(N+1)j \left( \int_{\lfloor \sqrt{j} \rfloor}^{\infty} \frac{1}{x^3} \, dx \right) \Delta t^2$$

$$\le C(N+1) \Delta t^2.$$

We have thus shown that $\sum_{k \in S_i^j} |E\Delta_k| \le C(N+1)\Delta t^2$.

Using (3.27)–(3.29), we have

$$E \left| \sum_{j=\lfloor ns \rfloor + 1}^{\lfloor nt \rfloor} \sigma_j^2 h_j (\sigma_j^{-1} \Delta F_j) \right|^4 \le C \sum_{j=0}^{N} (N+1) \Delta t^2 = C \left( \frac{\lfloor nt \rfloor - \lfloor ns \rfloor}{n} \right)^2,$$

which is (3.26).

To show that a sequence of cadlag processes $\{X_n\}$ is relatively compact, it suffices to show that for each $T > 1$, there exist constants $\beta > 0$, $C > 0$, and $\theta > 1$ such that

$$(3.30) \quad M_X(n, t, h) = E[|X_n(t+h) - X_n(t)|^\beta |X_n(t) - X_n(t-h)|^\beta] \leq Ch^\theta$$

for all $n \in \mathbb{N}$, all $t \in [0, T]$ and all $h \in [0, t]$. (See, e.g., Theorem 3.8.8 in [6].) Taking $\beta = 2$ and using (3.26) together with Hölder's inequality gives

$$M_B(n, t, h) \leq C\left(\frac{\lfloor nt + nh \rfloor - \lfloor nt \rfloor}{n}\right)\left(\frac{\lfloor nt \rfloor - \lfloor nt - nh \rfloor}{n}\right).$$

If $nh < 1/2$, then the right-hand side of this inequality is zero. Assume $nh \geq 1/2$. Then

$$\frac{\lfloor nt + nh \rfloor - \lfloor nt \rfloor}{n} \leq \frac{nh + 1}{n} \leq 3h.$$

The other factor is similarly bounded, so that $M_B(n, t, h) \leq Ch^2$.   □

Let us now introduce the filtration

$$\mathcal{F}_t = \sigma\{W(A) : A \subset \mathbb{R} \times [0, t], m(A) < \infty\},$$

where $m$ denotes Lebesgue measure on $\mathbb{R}^2$. Recall that

$$(3.31) \qquad F(t) = \int_{[0,t] \times \mathbb{R}} p(t - r, y) W(dr \times dy)$$

so that $F$ is adapted to $\{\mathcal{F}_t\}$. Also, given constants $0 \leq \tau \leq s \leq t$, we have

$$E[F(t) | \mathcal{F}_\tau] = \int_{[0,\tau] \times \mathbb{R}} p(t - r, y) W(dr \times dy)$$

and

$$E|E[F(t) - F(s) | \mathcal{F}_\tau]|^2 = \int_0^\tau \int_{\mathbb{R}} |p(t - r, y) - p(s - r, y)|^2 \, dy \, dr.$$

As in the proof of Lemma 2.1,

$$\int_0^\tau \int_{\mathbb{R}} p(t - r, y) p(s - r, y) \, dy \, dr = \frac{1}{\sqrt{2\pi}}(|t + s|^{1/2} - |(t - u) + (s - u)|^{1/2}).$$

Therefore, using (2.5), we can verify that

$$E|E[F(t) - F(s) | \mathcal{F}_\tau]|^2 = E|F(t) - F(s)|^2 - E|F(t - \tau) - F(s - \tau)|^2.$$

Combined with (2.3), this gives

$$E|E[F(t) - F(s) | \mathcal{F}_\tau]|^2 \leq \frac{2|t - s|^2}{|t - \tau|^{3/2}}.$$

In particular,

$$(3.32) \qquad E|E[\Delta F_j|\mathcal{F}_\tau]|^2 \leq \frac{2\Delta t^2}{(t_j - \tau)^{3/2}}$$

whenever $\tau \leq t_{j-1}$.

LEMMA 3.6. *Let $B_n$ be given by (3.25) and assume $\{h_j\}$ is independent of $\mathcal{F}_\infty$ and satisfies Assumption 3.1. Fix $0 \leq s < t$ and a constant $\kappa$. If*

$$\lim_{n\to\infty} E|B_n(t) - B_n(s)|^2 = \kappa^2(t-s),$$

*then*

$$B_n(t) - B_n(s) \Rightarrow \kappa|t-s|^{1/2}\chi$$

*as $n \to \infty$, where $\chi$ is a standard normal random variable.*

PROOF. We will prove the lemma by showing that every subsequence has a subsequence converging in law to the given random variable.

Let $\{n_j\}$ be any sequence. For each $n \in \mathbb{N}$, choose $m = m_n \in \{n_j\}$ such that $m_n > m_{n-1}$ and $m_n \geq n^4(t-s)^{-1}$. Now fix $n \in \mathbb{N}$ and let $\mu = m(t-s)/n$. For $0 \leq k < n$, define $u_k = ms + k\mu$, and let $u_n = mt$, so that

$$B_m(t) - B_m(s) = \sum_{j=\lfloor ms \rfloor + 1}^{mt} \sigma_j^2 h_j(\sigma_j^{-1}\Delta F_j)$$

$$= \sum_{k=1}^{n} \sum_{j=u_{k-1}+1}^{u_k} \sigma_j^2 h_j(\sigma_j^{-1}\Delta F_j).$$

For each pair $(j, k)$ such that $u_{k-1} < j \leq u_k$, let

$$\Delta\overline{F}_{j,k} = \Delta F_j - E[\Delta F_j|\mathcal{F}_{u_{k-1}\Delta t}].$$

Note that $\Delta\overline{F}_{j,k}$ is $\mathcal{F}_{u_k\Delta t}$-measurable and independent of $\mathcal{F}_{u_{k-1}\Delta t}$. We also make the following observation about $\Delta\overline{F}_{j,k}$. If we define

$$G_k(t) = F(t + \tau_k) - E[F(t + \tau_k)|\mathcal{F}_{\tau_k}],$$

where $\tau_k = u_{k-1}\Delta t$, then by (3.31),

$$G_k(t) = \int_{(\tau_k, t+\tau_k] \times \mathbb{R}} p(t + \tau_k - r, y)W(dr \times dy).$$

Hence, $G_k$ and $\mathcal{F}_{\tau_k}$ are independent, and $G_k$ and $F$ have the same law. Since

$$\Delta\overline{F}_{j,k} = \Delta F_j - E[\Delta F_j|\mathcal{F}_{\tau_k}] = G_k(t_j - \tau_k) - G_k(t_{j-1} - \tau_k),$$

it follows that $\{\Delta\overline{F}_{j,k}\}$ has the same law as $\{\Delta F_{j-u_{k-1}}\}$.

Now define $\bar{\sigma}_{j,k}^2 = E\Delta\overline{F}_{j,k}^2 = \sigma_{j-u_{k-1}}^2$ and

$$Z_{n,k} = \sum_{j=u_{k-1}+1}^{u_k} \bar{\sigma}_{j,k}^2 h_j(\bar{\sigma}_{j,k}^{-1}\Delta\overline{F}_{j,k})$$

so that $Z_{n,k}$, $1 \le k \le n$, are independent and

(3.33) $$B_m(t) - B_m(s) = \sum_{k=1}^n Z_{n,k} + \varepsilon_m,$$

where

$$\varepsilon_m = \sum_{k=1}^n \sum_{j=u_{k-1}+1}^{u_k} \{\sigma_j^2 h_j(\sigma_j^{-1}\Delta F_j) - \bar{\sigma}_{j,k}^2 h_j(\bar{\sigma}_{j,k}^{-1}\Delta\overline{F}_{j,k})\}.$$

Since $\Delta\overline{F}_{j,k}$ and $\Delta F_j - \Delta\overline{F}_{j,k} = E[\Delta F_j | \mathcal{F}_{u_{k-1}\Delta t}]$ are independent, we have

(3.34) $$\begin{aligned}\sigma_j^2 = E\Delta F_j^2 &= E\Delta\overline{F}_{j,k}^2 + E|\Delta F_j - \Delta\overline{F}_{j,k}|^2 \\ &= \bar{\sigma}_{j,k}^2 + E|\Delta F_j - \Delta\overline{F}_{j,k}|^2,\end{aligned}$$

which implies that $\bar{\sigma}_{j,k}^2 \le \sigma_j^2 \le C\Delta t^{1/2}$. In general, if $0 < a \le b$ and $x, y \in \mathbb{R}$, then by (3.5) and (3.9),

$$|b^2 h_j(b^{-1}y) - a^2 h_j(b^{-1}y)| \le (b^2 - a^2)L_j(1 + |b^{-1}y|^2),$$
$$|a^2 h_j(b^{-1}y) - a^2 h_j(a^{-1}x)| \le |a|^2 CL_j|b^{-1}y - a^{-1}x|(1 + |b^{-1}y| + |a^{-1}x|).$$

Note that $|b^{-1}y - a^{-1}x| \le |b^{-1} - a^{-1}||y| + |a^{-1}||y - x|$ and

$$|b^{-1} - a^{-1}| = \frac{b^2 - a^2}{ab(b+a)} \le \frac{b^2 - a^2}{a^3}.$$

Hence, if $\delta = b^2 - a^2$, then

$$\begin{aligned}|b^2 h_j&(b^{-1}y) - a^2 h_j(a^{-1}x)| \\ &\le CL_j(1 + |b^{-1}y|^2 + |a^{-1}x|)(\delta + \delta|a^{-1}y| + |a||y-x|).\end{aligned}$$

Using (2.9), Hölder's inequality and (3.34), this gives

$$\begin{aligned}E|\sigma_j^2 h_j&(\sigma_j^{-1}\Delta F_j) - \bar{\sigma}_{j,k}^2 h_j(\bar{\sigma}_{j,k}^{-1}\Delta\overline{F}_{j,k})| \\ &\le CE|\Delta F_j - \Delta\overline{F}_{j,k}|^2 + C\Delta t^{1/4}(E|\Delta F_j - \Delta\overline{F}_{j,k}|^2)^{1/2}.\end{aligned}$$

By (3.32),

$$E|\Delta F_j - \Delta\overline{F}_{j,k}|^2 \le \frac{2\Delta t^2}{(t_j - u_{k-1}\Delta t)^{3/2}} = \frac{2\Delta t^{1/2}}{(j - u_{k-1})^{3/2}}.$$

Therefore,

$$E|\varepsilon_m| \leq \sum_{k=1}^{n} \sum_{j=u_{k-1}+1}^{u_k} \frac{C\Delta t^{1/2}}{(j-u_{k-1})^{3/4}} = Cm^{-1/2} \sum_{k=1}^{n} \sum_{j=1}^{u_k-u_{k-1}} j^{-3/4}.$$

Since $u_k - u_{k-1} \leq C\mu$, this gives

$$E|\varepsilon_m| \leq Cm^{-1/2} n\mu^{1/4} = Cn^{3/4}m^{-1/4}(t-s)^{1/4}.$$

But since $m = m_n$ was chosen so that $m \geq n^4(t-s)^{-1}$, we have $E|\varepsilon_m| \leq Cn^{-1/4}|t-s|^{1/2}$ and $\varepsilon_m \to 0$ in $L^1$ and in probability. Therefore, by (3.33), we need only to show that

$$\sum_{k=1}^{n} Z_{n,k} \Rightarrow \kappa|t-s|^{1/2}\chi$$

in order to complete the proof.

For this, we will use the Lindeberg–Feller theorem (see, e.g. Theorem 2.4.5 in [5]), which states the following: for each $n$, let $Z_{n,k}$, $1 \leq k \leq n$, be independent random variables with $EZ_{n,k} = 0$. Suppose:

(a) $\sum_{k=1}^{n} EZ_{n,k}^2 \to \sigma^2$, and

(b) for all $\varepsilon > 0$, $\lim_{n\to\infty} \sum_{k=1}^{n} E[|Z_{n,k}|^2 1_{\{|Z_{n,k}|>\varepsilon\}}] = 0$.

Then $\sum_{k=1}^{n} Z_{n,k} \Rightarrow \sigma\chi$ as $n \to \infty$.

To verify these conditions, recall that $\{\Delta\overline{F}_{j,k}\}$ and $\{\Delta F_{j-u_{k-1}}\}$ have the same law, so that

$$E|Z_{n,k}|^4 = E\left| \sum_{j=u_{k-1}+1}^{u_k} \bar\sigma_{j,k}^2 h_j(\bar\sigma_{j,k}^{-1}\Delta\overline{F}_{j,k}) \right|^4$$

$$= E\left| \sum_{j=1}^{u_k-u_{k-1}} \sigma_j^2 h_{j+u_{k-1}}(\sigma_j^{-1}\Delta F_j) \right|^4$$

$$= E|\bar{B}_{m,k}((u_k-u_{k-1})\Delta t)|^4,$$

where

$$\bar{B}_{m,k}(t) = \sum_{j=1}^{\lfloor mt \rfloor} \sigma_j^2 h_{j+u_{k-1}}(\sigma_j^{-1}\Delta F_j).$$

Hence, by Proposition 3.5,

$$E|Z_{n,k}|^4 \leq C(u_k - u_{k-1})^2 \Delta t^2.$$

Jensen's inequality now gives $\sum_{k=1}^{n} E|Z_{n,k}|^2 \leq Cn\mu\Delta t = C(t-s)$, so that by passing to a subsequence, we may assume that (a) holds for some $\sigma \geq 0$.

For (b), let $\varepsilon > 0$ be arbitrary. Then

$$\sum_{k=1}^{n} E[|Z_{n,k}|^2 1_{\{|Z_{n,k}|>\varepsilon\}}] \le \varepsilon^{-2} \sum_{k=1}^{n} E|Z_{n,k}|^4$$
$$\le C\varepsilon^{-2} n\mu^2 \Delta t^2$$
$$= C\varepsilon^{-2} n^{-1}(t-s)^2,$$

which tends to zero as $n \to \infty$.

It therefore follows that $\sum_{k=1}^{n} Z_{n,k} \Rightarrow \sigma\chi$ as $n \to \infty$ and it remains only to show that $\sigma = \kappa|t-s|^{1/2}$. For this, observe that the continuous mapping theorem implies that $|B_m(t) - B_m(s)|^2 \Rightarrow \sigma^2 \chi^2$. By the Skorohod representation theorem, we may assume that the convergence is a.s. By Proposition 3.5, the family $|B_m(t) - B_m(s)|^2$ is uniformly integrable. Hence, $|B_m(t) - B_m(s)|^2 \to \sigma^2 \chi^2$ in $L^1$, which implies $E|B_m(t) - B_m(s)|^2 \to \sigma^2$. But by assumption, $E|B_m(t) - B_m(s)|^2 \to \kappa^2(t-s)$, so $\sigma = \kappa|t-s|^{1/2}$ and the proof is complete. $\square$

LEMMA 3.7.  *Let $B_n$ be given by (3.25) and assume $\{h_j\}$ is independent of $\mathcal{F}_\infty$ and satisfies Assumption 3.1, so that by Proposition 3.5, the sequence $\{B_n\}$ is relatively compact. If $X$ is any weak limit point of this sequence, then $X$ has independent increments.*

PROOF.  Suppose that $B_{n(j)} \Rightarrow X$. Fix $0 < t_1 < t_2 < \cdots < t_d < s < t$. It will be shown that $X(t) - X(s)$ and $(X(t_1), \ldots, X(t_d))$ are independent. With notation as in Lemma 3.6, let

$$Z_n = \sum_{j=\lfloor ns \rfloor+2}^{\lfloor nt \rfloor} \bar{\sigma}_{j,k}^2 h_j(\bar{\sigma}_{j,k}^{-1} \Delta \overline{F}_{j,k}),$$

and define

$$Y_n = B_n(t) - B_n(s) - Z_n.$$

As in the proof of Lemma 3.6, $Y_n \to 0$ in probability. It therefore follows that

$$(B_{n(j)}(t_1), \ldots, B_{n(j)}(t_d), Z_{n(j)}) \Rightarrow (X(t_1), \ldots, X(t_d), X(t) - X(s)).$$

Note that $\mathcal{F}_{(\lfloor ns \rfloor+1)\Delta t}$ and $Z_n$ are independent. Hence, $(B_n(t_1), \ldots, B_n(t_d))$ and $Z_n$ are independent, which implies $X(t) - X(s)$ and $(X(t_1), \ldots, X(t_d))$ are independent.  $\square$

THEOREM 3.8.  *Let*

$$(3.35) \qquad B_n(t) = \sum_{j=1}^{\lfloor nt \rfloor} \sigma_j^2 h_j(\sigma_j^{-1} \Delta F_j)$$

*and assume $\{h_j\}$ is independent of $\mathcal{F}_\infty$ and satisfies Assumption 3.1. If there exists a constant $\kappa$ such that*

$$\lim_{n\to\infty} E|B_n(t) - B_n(s)|^2 = \kappa^2(t-s)$$

*for all $0 \le s < t$, then $(F, B_n) \Rightarrow (F, \kappa B)$, where $B$ is a standard Brownian motion independent of $F$.*

PROOF. Let $\{n(j)\}_{j=1}^\infty$ be any sequence of natural numbers. By Proposition 3.5, the sequence $\{(F, B_{n(j)})\}$ is relatively compact. Therefore, there exists a subsequence $m(k) = n(j_k)$ and a cadlag process $X$ such that $(F, B_{m(k)}) \Rightarrow (F, X)$. By Lemma 3.7, the process $X$ has independent increments. By Lemma 3.6, the increment $X(t) - X(s)$ is normally distributed with mean zero and variance $\kappa^2|t-s|$. Also, $X(0) = 0$ since $B_n(0) = 0$ for all $n$. Hence, $X$ is equal in law to $\kappa B$, where $B$ is a standard Brownian motion. It remains only to show that $F$ and $B$ are independent.

Fix $0 < r_1 < \cdots < r_\ell$ and define $\xi = (F(r_1), \ldots, F(r_\ell))^T$. It is easy to see that $\Sigma = E\xi\xi^T$ is invertible. Hence, we may define the vectors $v_j \in \mathbb{R}^\ell$ by $v_j = E[\xi \Delta F_j]$, and $a_j = \Sigma^{-1} v_j$. Let $\Delta^* F_j = \Delta F_j - a_j^T \xi$, so that $\Delta^* F_j$ and $\xi$ are independent.

Define

$$B_n^*(t) = \sum_{j=1}^{\lfloor nt \rfloor} \sigma_j^2 h_j(\sigma_j^{-1} \Delta^* F_j).$$

As in the proof of Lemma 3.6,

$$E\left[\sup_{0 \le t \le T} |B_n(t) - B_n^*(t)|\right] \le C \sum_{j=1}^{\lfloor nT \rfloor} |a_j| \Delta t^{1/4} + C \sum_{j=1}^{\lfloor nT \rfloor} |a_j|^2,$$

where $C$ is a constant that depends only on $(r_1, \ldots, r_\ell)$. Also note that

$$|a_j| \le C|v_j| \le C \sum_{i=1}^\ell |E[F(r_i)\Delta F_j]|.$$

Hence,

$$E\left[\sup_{0 \le t \le T} |B_n(t) - B_n^*(t)|\right] \le C\Delta t^{1/4} \sum_{i=1}^\ell \sum_{j=1}^{\lfloor nT \rfloor} |E[F(r_i)\Delta F_j]|$$

$$+ C \sum_{i=1}^\ell \sum_{j=1}^{\lfloor nT \rfloor} |E[F(r_i)\Delta F_j]|^2.$$

Using estimates similar to those in the proof of (2.11), it can be verified that

$$\limsup_{n \to \infty} \sum_{j=1}^{\lfloor nt \rfloor} |E[F(r_i)\Delta F_j]| \le C \int_0^t \frac{1}{|u - r_i|^{1/2}} \, du < \infty$$

and

$$\lim_{n \to \infty} \sum_{j=1}^{\lfloor nt \rfloor} |E[F(r_i)\Delta F_j]|^2 = 0.$$

Thus, $(\xi, B_n^*(r_1), \dots, B_n^*(r_\ell)) \Rightarrow (\xi, \kappa B(r_1), \dots, \kappa B(r_\ell))$. Since $\xi$ and $B_n^*$ are independent, this finishes the proof. $\quad\square$

## 4. Examples.

4.1. *Independent mean zero sign changes.*

PROPOSITION 4.1. *Let $\{\xi_j\}$ be a sequence of independent mean zero random variables with $E\xi_j^2 = 1$. Suppose that the sequence $\{\xi_j\}$ is independent of $\mathcal{F}_\infty$. Let*

$$B_n(t) = \sum_{j=1}^{\lfloor nt \rfloor} \Delta F_j^2 \xi_j.$$

*Then $(F, B_n) \Rightarrow (F, 6\pi^{-1}B)$, where $B$ is a standard Brownian motion independent of $F$.*

PROOF. Let $h_j(x) = \xi_j x^2$. Then $\{h_j\}$ satisfies Assumption 3.1 with $L_j = 2|\xi_j|$ and $L = 4$, and $B_n$ has the form (3.35). Moreover,

$$E|B_n(t) - B_n(s)|^2 = E\left| \sum_{j=\lfloor ns \rfloor + 1}^{\lfloor nt \rfloor} \Delta F_j^2 \xi_j \right|^2 = \sum_{j=\lfloor ns \rfloor + 1}^{\lfloor nt \rfloor} E\Delta F_j^4.$$

By (2.12),

$$\lim_{n \to \infty} E|B_n(t) - B_n(s)|^2 = 6\pi^{-1}(t - s).$$

The result now follows from Theorem 3.8. $\quad\square$

4.2. *The signed variations of $F$.* In this subsection, we adopt the notation $x^{r\pm} = |x|^r \operatorname{sgn}(r)$. We begin by showing that the "signed cubic variation" of $F$ is zero.

PROPOSITION 4.2. *If $Z_n(t) = \sum_{j=1}^{\lfloor nt \rfloor} \Delta F_j^3$, then $Z_n(t) \to 0$ uniformly on compacts in probability.*

PROOF. Note that $x_n \to 0$ in $D_{\mathbb{R}}[0, \infty)$ if and only if $x_n \to 0$ uniformly on compacts. Hence, we must show that $Z_n \to 0$ in probability in $D_{\mathbb{R}}[0, \infty)$, for which it will suffice to show that $Z_n \Rightarrow 0$.

Note that

$$E|Z_n(t) - Z_n(s)|^2 = E \left| \sum_{j=\lfloor ns \rfloor + 1}^{\lfloor nt \rfloor} \Delta F_j^3 \right|^2$$

$$\leq \sum_{i=\lfloor ns \rfloor + 1}^{\lfloor nt \rfloor} \sum_{j=\lfloor ns \rfloor + 1}^{\lfloor nt \rfloor} |E[\Delta F_i^3 \Delta F_j^3]|.$$

To estimate this sum, we use the following fact about Gaussian random variables. Let $X_1, X_2$ be mean zero, jointly normal random variables with variances $\sigma_j^2$. If $\rho = (\sigma_1 \sigma_2)^{-1} E[X_1 X_2]$, then

$$E[X_1^3 X_2^3] = \sigma_1^3 \sigma_2^3 \rho (6\rho^2 + 9).$$

Applying this in our context, let $\rho_{ij} = (\sigma_i \sigma_j)^{-1} E[\Delta F_i \Delta F_j]$, so that

$$E|Z_n(t) - Z_n(s)|^2 \leq C \sum_{i=ns+1}^{\lfloor nt \rfloor} \sum_{j=\lfloor ns \rfloor + 1}^{\lfloor nt \rfloor} \sigma_i^3 \sigma_j^3 |\rho_{ij}|$$

$$= C \sum_{i=ns+1}^{\lfloor nt \rfloor} \sum_{j=\lfloor ns \rfloor + 1}^{\lfloor nt \rfloor} \sigma_i^2 \sigma_j^2 |E[\Delta F_i \Delta F_j]|.$$

Using (2.10), this gives

$$E|Z_n(t) - Z_n(s)|^2 \leq C \sum_{i=\lfloor ns \rfloor + 1}^{\lfloor nt \rfloor} \sum_{j=\lfloor ns \rfloor + 1}^{\lfloor nt \rfloor} \sqrt{\Delta t} \sqrt{\Delta t} \left( \frac{\sqrt{\Delta t}}{|i-j|^{\sim 3/2}} \right)$$

$$\leq C \left( \frac{\lfloor nt \rfloor - \lfloor ns \rfloor}{n} \right) \sqrt{\Delta t}.$$

Hence, $Z_n(t) \to 0$ in probability for each fixed $t$. Moreover, taking $\beta = 1$ in (3.30), this shows that $M_Z(n, t, h) = 0$ when $nh < 1/2$, and $M_Z(n, t, h) \leq Ch\sqrt{\Delta t} \leq Ch^{3/2}$ when $nh \geq 1/2$. Therefore, $\{Z_n\}$ is relatively compact and $Z_n \Rightarrow 0$. $\quad\square$

LEMMA 4.3. *Let $X_1, X_2$ be mean zero, jointly normal random variables with $EX_j^2 = 1$ and $\rho = E[X_1 X_2]$. Let*

$$K(x) = \frac{6}{\pi} x \sqrt{1 - x^2} + \frac{2}{\pi} (1 + 2x^2) \sin^{-1}(x),$$

*where $\sin^{-1}(x) \in [-\pi/2, \pi/2]$. Then $E[X_1^{2\pm} X_2^{2\pm}] = K(\rho)$. Moreover, for all $x \in [-1, 1]$, we have $|K(x) - 8x/\pi| \leq 2|x|^3$, so that $|E[X_1^{2\pm} X_2^{2\pm}]| \leq 5|\rho|$.*

Proof. Define $U = X_1$ and $V = (1 - \rho^2)^{-1/2}(X_2 - \rho X_1)$, so that $U$ and $V$ are independent standard normals. Then $X_1 = U$ and $X_2 = \eta V + \rho U$, where $\eta = \sqrt{1 - \rho^2}$, and

$$
\begin{aligned}
E[X_1^{2\pm} X_2^{2\pm}] &= \frac{1}{2\pi} \iint [u(\eta v + \rho u)]^{2\pm} e^{-(u^2 + v^2)/2} \, du \, dv \\
&= \frac{1}{2\pi} \int_0^{2\pi} \int_0^\infty [\cos\theta(\eta \sin\theta + \rho \cos\theta)]^{2\pm} r^5 e^{-r^2/2} \, dr \, d\theta \\
&= \frac{4}{\pi} \int_0^{2\pi} [\cos^2\theta(\eta \tan\theta + \rho)]^{2\pm} \, d\theta.
\end{aligned}
$$

If $a = \tan^{-1}(-\rho/\eta)$, then we can write

$$
\begin{aligned}
E[X_1^{2\pm} X_2^{2\pm}] &= \frac{8}{\pi} \int_a^{\pi/2} [\cos^2\theta(\eta \tan\theta + \rho)]^2 \, d\theta \\
&\quad - \frac{8}{\pi} \int_{-\pi/2}^a [\cos^2\theta(\eta \tan\theta + \rho)]^2 \, d\theta \\
&= \frac{8}{\pi} \int_{-\pi/2}^{-a} [\cos^2\theta(\eta \tan\theta - \rho)]^2 \, d\theta \\
&\quad - \frac{8}{\pi} \int_{-\pi/2}^a [\cos^2\theta(\eta \tan\theta + \rho)]^2 \, d\theta.
\end{aligned}
$$

By symmetry, we can assume that $\rho \leq 0$, so that $a \geq 0$. Then

$$
\begin{aligned}
E[X_1^{2\pm} X_2^{2\pm}] &= -\frac{32}{\pi} \rho\eta \int_{-\pi/2}^{-a} \cos^4\theta \tan\theta \, d\theta - \frac{8}{\pi} \int_{-a}^a [\cos^2\theta(\eta \tan\theta + \rho)]^2 \, d\theta \\
&= -\frac{32}{\pi} \rho\eta \left(-\frac{1}{4} \cos^4 a\right) - \frac{16}{\pi} \int_0^a \cos^4\theta \, d\theta \\
&\quad + \frac{16}{\pi} \eta^2 \int_0^a \cos^4\theta(1 - \tan^2\theta) \, d\theta.
\end{aligned}
$$

Using $a = \sin^{-1}(-\rho)$ and the formulas

$$
\int \cos^4\theta(1 - \tan^2\theta) \, d\theta = (\theta + \sin\theta\cos\theta + 2\sin\theta\cos^3\theta)/4,
$$

$$
\int \cos^4\theta \, d\theta = (3\theta + 3\sin\theta\cos\theta + 2\sin\theta\cos^3\theta)/8,
$$

we can directly verify that $E[X_1^{2\pm} X_2^{2\pm}] = K(\rho)$.

To estimate $K$, note that $K \in C^\infty(-1, 1)$ with

$$
K'(x) = \frac{8}{\pi}(\sqrt{1 - x^2} + x \sin^{-1}(x)),
$$

$$
K''(x) = \frac{8}{\pi} \sin^{-1}(x).
$$

Since $K''$ is increasing,

$$\left| K(x) - \frac{8}{\pi}x \right| \leq \frac{1}{2}x^2 K''(|x|).$$

But for $y \in [0, \pi/2]$, we have $\sin y \geq 2y/\pi$. Letting $y = \pi x/2$ gives $\sin^{-1}(x) \leq \pi x/2$ for $x \in [0, 1]$. We therefore have $K''(|x|) \leq 4|x|$, so that $|K(x) - 8x/\pi| \leq 2|x|^3$. $\square$

PROPOSITION 4.4. *Let $K$ be defined as in Lemma 4.3, and $\gamma_i$ as in Lemma 2.1. If*

$$B_n(t) = \sum_{j=1}^{\lfloor nt \rfloor} \Delta F_j^2 \, \mathrm{sgn}(\Delta F_j),$$

*then $(F, B_n) \Rightarrow (F, \kappa B)$, where $\kappa^2 = 6\pi^{-1} - 4\pi^{-1} \sum_{i=0}^{\infty} K(\gamma_i/2) > 0$ and $B$ is a standard Brownian motion independent of $F$.*

PROOF. Let $h_j(x) = h(x) = x^{2\pm}$, so that

$$B_n(t) = \sum_{j=1}^{\lfloor nt \rfloor} \sigma_j^2 h_j(\sigma_j^{-1} \Delta F_j).$$

Since $h$ is continuously differentiable and $h'(x) = 2|x|$ is Lipschitz, $\{h_j\}$ satisfies (3.1). Moreover, if $X$ and $Y$ are jointly normal with mean zero, variance one, and covariance $\rho = EXY$, then $Eh(X) = 0$ and $|Eh(X)h(Y)| \leq 5|\rho|$ by Lemma 4.3. Hence, $\{h_j\}$ satisfies Assumption 3.1. By Proposition 3.5,

$$(4.1) \qquad E|B_n(t) - B_n(s)|^4 \leq C\left( \frac{\lfloor nt \rfloor - \lfloor ns \rfloor}{n} \right)^2$$

for all $0 \leq s < t$ and all $n \in \mathbb{N}$.

By Theorem 3.8, the proof will be complete once we establish that $\kappa$ is well-defined, strictly positive and

$$(4.2) \qquad \lim_{n \to \infty} E|B_n(t) - B_n(s)|^2 = \kappa^2(t - s)$$

for all $0 \leq s < t$.

By (2.8), $\gamma_i/2 \in (0, 1]$ for all $i$. Thus, by Lemma 4.3,

$$0 < \sum_{i=1}^{\infty} K(\gamma_i/2) \leq \frac{8}{\pi} \sum_{i=1}^{\infty} \frac{1}{2}\gamma_i + 2\sum_{i=1}^{\infty} \left( \frac{1}{2}\gamma_i \right)^3 = \frac{4}{\pi} \sum_{i=1}^{\infty} \gamma_i + \frac{1}{4} \sum_{i=1}^{\infty} \gamma_i^3.$$

Since $\gamma_i = f(i-1) - f(i)$, where $f(x) = \sqrt{x+1} - \sqrt{x}$, we have that $\sum_{i=1}^{\infty} \gamma_i = f(0) = 1$. Moreover, by (2.8),

$$\sum_{i=1}^{\infty} \gamma_i^3 \leq \sum_{i=1}^{\infty} \frac{1}{(\sqrt{2}\, i^{3/2})^3} = \frac{1}{2\sqrt{2}} \sum_{i=1}^{\infty} \frac{1}{i^{9/2}}$$

$$\leq \frac{1}{2\sqrt{2}}\left(1 + \int_1^\infty \frac{1}{x^{9/2}}\,dx\right) = \frac{9\sqrt{2}}{28}.$$

Thus,

$$\sum_{i=1}^\infty K(\gamma_i/2) \leq \frac{4}{\pi} + \frac{9\sqrt{2}}{112} < \frac{3}{2},$$

which gives $6\pi^{-1} - 4\pi^{-1}\sum_{i=0}^\infty K(\gamma_i/2) > 0$, so that $\kappa$ is well-defined and strictly positive.

Now fix $0 \leq s < t$. First assume that $s > 0$. Then

$$E|B_n(t) - B_n(s)|^2 = E\left|\sum_{j=\lfloor ns\rfloor+1}^{\lfloor nt\rfloor} \Delta F_j^{2\pm}\right|^2$$

$$= \sum_{j=\lfloor ns\rfloor+1}^{\lfloor nt\rfloor} E\Delta F_j^4 + 2\sum_{j=\lfloor ns\rfloor+2}^{\lfloor nt\rfloor}\sum_{i=\lfloor ns\rfloor+1}^{j-1} E[\Delta F_i^{2\pm}\Delta F_j^{2\pm}]$$

$$= \sum_{j=\lfloor ns\rfloor+1}^{\lfloor nt\rfloor} \frac{6}{\pi}\Delta t - \frac{4}{\pi}\sum_{j=\lfloor ns\rfloor+2}^{\lfloor nt\rfloor}\sum_{i=\lfloor ns\rfloor+1}^{j-1} K(\gamma_{j-i}/2)\Delta t + R_n,$$

where

$$R_n = \sum_{j=\lfloor ns\rfloor+1}^{\lfloor nt\rfloor}\left(E\Delta F_j^4 - \frac{6}{\pi}\Delta t\right)$$

$$+ 2\sum_{j=\lfloor ns\rfloor+2}^{\lfloor nt\rfloor}\sum_{i=\lfloor ns\rfloor+1}^{j-1}\left(E[\Delta F_i^{2\pm}\Delta F_j^{2\pm}] + \frac{2}{\pi}K(\gamma_{j-i}/2)\Delta t\right).$$

Observe that

$$\sum_{j=\lfloor ns\rfloor+1}^{\lfloor nt\rfloor} \frac{6}{\pi}\Delta t = \frac{6}{\pi}\left(\frac{\lfloor nt\rfloor - \lfloor ns\rfloor}{n}\right) \to \frac{6}{\pi}(t-s)$$

and

$$\sum_{j=\lfloor ns\rfloor+2}^{\lfloor nt\rfloor}\sum_{i=\lfloor ns\rfloor+1}^{j-1} K(\gamma_{j-i}/2)\Delta t = \sum_{j=\lfloor ns\rfloor+2}^{\lfloor nt\rfloor}\sum_{i=1}^{j-\lfloor ns\rfloor-1} K(\gamma_i/2)\Delta t$$

(4.3)

$$= \sum_{j=1}^N \sum_{i=1}^j \frac{1}{n} K(\gamma_i/2),$$

where $N = \lfloor nt \rfloor - \lfloor ns \rfloor - 1$. Thus,

$$
\sum_{j=\lfloor ns \rfloor + 2}^{\lfloor nt \rfloor} \sum_{i=\lfloor ns \rfloor + 1}^{j-1} K(\gamma_{j-i}/2) \Delta t = \sum_{i=1}^{N} \sum_{j=i}^{N} \frac{1}{n} K(\gamma_i/2)
$$

$$(4.4)$$

$$
= \sum_{i=1}^{N} \left( \frac{N}{n} - \frac{i}{n} \right) K(\gamma_i/2).
$$

We claim that $n^{-1} \sum_{i=1}^{N} i K(\gamma_i/2) \to 0$ as $n \to \infty$. To see this, fix $m \in \mathbb{N}$ and note that

$$
\frac{1}{n} \sum_{i=1}^{N} i K(\gamma_i/2) \leq \frac{m}{n} \sum_{i=1}^{m} K(\gamma_i/2) + \frac{N}{n} \sum_{i=m+1}^{N} K(\gamma_i/2).
$$

Since $N/n \to (t-s)$ as $n \to \infty$, this gives

$$
\limsup_{n \to \infty} \frac{1}{n} \sum_{i=1}^{N} i K(\gamma_i/2) \leq (t-s) \sum_{i=m+1}^{\infty} K(\gamma_i/2).
$$

Letting $m \to \infty$ proves the claim. It now follows that

$$
\sum_{j=\lfloor ns \rfloor + 1}^{\lfloor nt \rfloor} \frac{6}{\pi} \Delta t - \frac{4}{\pi} \sum_{j=\lfloor ns \rfloor + 2}^{\lfloor nt \rfloor} \sum_{i=\lfloor ns \rfloor + 1}^{j-1} K(\gamma_{j-i}/2) \Delta t \to \kappa^2(t-s)
$$

and it suffices to show that $R_n \to 0$ as $n \to \infty$.

Now, by Lemma 4.3,

$$
E[\Delta F_i^{2\pm} \Delta F_j^{2\pm}] = \sigma_i^2 \sigma_j^2 E\left[ \left( \frac{\Delta F_i}{\sigma_i} \right)^{2\pm} \left( \frac{\Delta F_j}{\sigma_j} \right)^{2\pm} \right] = \sigma_i^2 \sigma_j^2 K(\rho_{ij})
$$

where $\sigma_i^2 = E\Delta F_i^2$ and

$$
\rho_{ij} = E\left[ \left( \frac{\Delta F_i}{\sigma_i} \right) \left( \frac{\Delta F_j}{\sigma_j} \right) \right] = (\sigma_i \sigma_j)^{-1} E[\Delta F_i \Delta F_j].
$$

Define $a = \sigma_i^2 \sigma_j^2$, $b = K(\rho_{ij})$, $c = 2\Delta t/\pi$, and $d = K(-\gamma_{j-i}/2)$. By (2.10), we have $|a| \leq C\Delta t$. By (2.10) and (2.3),

$$
|a - c| = \left| \sigma_i^2 \left( \sigma_j^2 - \sqrt{\frac{2\Delta t}{\pi}} \right) + \sqrt{\frac{2\Delta t}{\pi}} \left( \sigma_i^2 - \sqrt{\frac{2\Delta t}{\pi}} \right) \right| \leq C \frac{1}{t_i^{3/2}} \Delta t^{5/2}.
$$

By Lemma 4.3, $|K(x) - K(y)| \leq C|x-y|$, so that $|d| \leq C$ and $|b-d| \leq C|\rho_{ij} + \gamma_{j-i}/2|$. Rewriting this latter inequality, we have that $|b-d|$ is bounded above by

$$
C\left| \frac{1}{\sigma_i \sigma_j} \left( E[\Delta F_i \Delta F_j] + \sqrt{\frac{\Delta t}{2\pi}} \gamma_{j-i} \right) - \sqrt{\frac{\Delta t}{2\pi}} \gamma_{j-i} \left( \frac{1}{\sigma_i \sigma_j} - \sqrt{\frac{\pi}{2\Delta t}} \right) \right|.
$$

Observe that

$$\sigma_i \sigma_j \left( \frac{1}{\sigma_i \sigma_j} - \sqrt{\frac{\pi}{2\Delta t}} \right)$$

$$= \sqrt{\frac{\pi}{2\Delta t}} \left[ \frac{\sigma_j}{\sigma_i + \sigma_j} \left( \sqrt{\frac{2\Delta t}{\pi}} - \sigma_i^2 \right) + \frac{\sigma_i}{\sigma_i + \sigma_j} \left( \sqrt{\frac{2\Delta t}{\pi}} - \sigma_j^2 \right) \right]$$

so that by (2.10) and (2.3)

$$\left| \sigma_i^2 \sigma_j^2 \left( \frac{1}{\sigma_i \sigma_j} - \sqrt{\frac{\pi}{2\Delta t}} \right) \right| \le C \frac{1}{t_i^{3/2}} \Delta t^2.$$

Hence, by (2.4), $|a||b - d| \le C t_i^{-3/2} \Delta t^{5/2}$. We therefore have

$$\left| E[\Delta F_i^{2\pm} \Delta F_j^{2\pm}] + \frac{2}{\pi} K(\gamma_{j-i}/2) \Delta t \right| = |ab - cd|$$

$$\le |a||b - d| + |d||a - c|$$

$$\le C t_i^{-3/2} \Delta t^{5/2}.$$

Since $t_i > s > 0$, this shows that

$$\sum_{j=\lfloor ns \rfloor + 2}^{\lfloor nt \rfloor} \sum_{i=\lfloor ns \rfloor + 1}^{j-1} \left( E[\Delta F_i^{2\pm} \Delta F_j^{2\pm}] + \frac{2}{\pi} K(\gamma_{j-i}/2) \Delta t \right) \to 0.$$

Combined with (2.12), this shows that $R_n \to 0$.

We have now proved (4.2) under the assumption that $s > 0$. Now assume $s = 0$. Let $\varepsilon \in (0, t)$ be arbitrary. Then by Hölder's inequality and (4.1),

$$|E|B_n(t)|^2 - \kappa^2 t| = |E|B_n(t) - B_n(\varepsilon)|^2 - \kappa^2(t - \varepsilon)$$

$$+ 2E[B_n(t)B_n(\varepsilon)] - E|B_n(\varepsilon)|^2 - \kappa^2 \varepsilon|$$

$$\le |E|B_n(t) - B_n(\varepsilon)|^2 - \kappa^2(t - \varepsilon)| + C(\sqrt{t\varepsilon} + \varepsilon).$$

First let $n \to \infty$, then let $\varepsilon \to 0$, and the proof is complete.  $\square$

### 4.3. *Centering the squared increments.*

PROPOSITION 4.5.   *Let* $\gamma_i$ *be defined as in Lemma 2.1. If*

$$B_n(t) = \sum_{j=1}^{\lfloor nt \rfloor} (\Delta F_j^2 - \sigma_j^2),$$

*then* $(F, B_n) \Rightarrow (F, \kappa B)$*, where* $\kappa^2 = 4\pi^{-1} + 2\pi^{-1} \sum_{i=0}^{\infty} \gamma_i^2$ *and* $B$ *is a standard Brownian motion independent of* $F$.

PROOF.   Let $h_j(x) = x^2 - 1$. Then $\{h_j\}$ clearly satisfies (3.1) and (3.2). For jointly normal $X$ and $Y$ with mean zero and variance one, $E(X^2 - 1)(Y^2 - 1) = 2\rho^2$, so $\{h_j\}$ also satisfies (3.3). Since

$$B_n(t) = \sum_{j=1}^{\lfloor nt \rfloor} \sigma_j^2 h_j(\sigma_j^{-1} \Delta F_j),$$

it will suffice, by Theorem 3.8, to show that

$$(4.5) \qquad \lim_{n \to \infty} E|B_n(t) - B_n(s)|^2 = \kappa^2(t - s).$$

By Proposition 3.5,

$$E|B_n(t) - B_n(s)|^4 \leq C\left(\frac{\lfloor nt \rfloor - \lfloor ns \rfloor}{n}\right)^2$$

for all $0 \leq s < t$ and all $n \in \mathbb{N}$. Hence, as in the proof of Proposition 4.4, it will suffice to prove (4.5) for $s > 0$.

Assume $s > 0$. Then

$$E|B_n(t) - B_n(s)|^2 = E\left|\sum_{j=\lfloor ns \rfloor + 1}^{\lfloor nt \rfloor} (\Delta F_j^2 - \sigma_j^2)\right|^2$$

$$= \sum_{j=\lfloor ns \rfloor + 1}^{\lfloor nt \rfloor} 2\sigma_j^4 + 2 \sum_{j=\lfloor ns \rfloor + 2}^{\lfloor nt \rfloor} \sum_{i=\lfloor ns \rfloor + 1}^{j-1} 2|E\Delta F_i \Delta F_j|^2$$

$$= \frac{4}{\pi}\left(\frac{\lfloor nt \rfloor - \lfloor ns \rfloor}{n}\right) + \frac{2}{\pi} \sum_{j=\lfloor ns \rfloor + 2}^{\lfloor nt \rfloor} \sum_{i=\lfloor ns \rfloor + 1}^{j-1} \gamma_{j-i}^2 \Delta t + R_n,$$

where

$$R_n = \sum_{j=\lfloor ns \rfloor + 1}^{\lfloor nt \rfloor} \left(2\sigma_j^4 - \frac{4}{\pi}\Delta t\right)$$

$$+ 4 \sum_{j=\lfloor ns \rfloor + 2}^{\lfloor nt \rfloor} \sum_{i=\lfloor ns \rfloor + 1}^{j-1} \left(|E\Delta F_i \Delta F_j|^2 - \frac{\Delta t}{2\pi}\gamma_{j-i}^2\right).$$

By (2.4), (2.8) and (2.10),

$$\left||E\Delta F_i \Delta F_j|^2 - \frac{\Delta t}{2\pi}\gamma_{j-i}^2\right| \leq \frac{C\Delta t^{5/2}}{s^{3/2}(j-i)^{3/2}}.$$

As in (4.3) and (4.4),

$$\sum_{j=\lfloor ns \rfloor + 2}^{\lfloor nt \rfloor} \sum_{i=\lfloor ns \rfloor + 1}^{j-1} \frac{\Delta t}{(j-i)^{3/2}} \to (t - s) \sum_{i=1}^{\infty} \frac{1}{i^{3/2}}.$$

Together with (2.12), this shows that $R_n \to 0$. Hence,

$$E|B_n(t) - B_n(s)|^2 \to \frac{4}{\pi}(t-s) + \frac{2}{\pi}(t-s)\sum_{i=1}^{\infty} \gamma_i^2,$$

and the proof is complete.  □

Corollary 4.6.   *If*

$$B_n(t) = \left(\sum_{j=1}^{\lfloor nt \rfloor} \Delta F_j^2\right) - \sqrt{\frac{2n}{\pi}}\, t$$

*then* $(F, B_n) \Rightarrow (F, \kappa B)$, *where* $\kappa^2 = 4\pi^{-1} + 2\pi^{-1}\sum_{i=0}^{\infty} \gamma_i^2$ *and* $B$ *is a standard Brownian motion independent of* $F$.

Proof.   Note that

$$B_n(t) = \sqrt{\frac{2}{\pi n}}(\lfloor nt \rfloor - nt) + \sum_{j=1}^{\lfloor nt \rfloor}\left(\Delta F_j^2 - \sqrt{\frac{2\Delta t}{\pi}}\right),$$

and by (2.3),

$$\sup_{0 \le s \le t}\left|\sum_{j=1}^{\lfloor ns \rfloor}\left(\Delta F_j^2 - \sqrt{\frac{2\Delta t}{\pi}}\right) - \sum_{j=1}^{\lfloor ns \rfloor}(\Delta F_j^2 - \sigma_j^2)\right| \le \sum_{j=1}^{\lfloor nt \rfloor}\left|\sigma_j^2 - \sqrt{\frac{2\Delta t}{\pi}}\right|$$

$$\le \sqrt{\Delta t}\sum_{j=1}^{\lfloor nt \rfloor}\frac{1}{j^{3/2}}.$$

The result now follows from Proposition 4.5.  □

### 4.4.  *Alternating sign changes.*

Proposition 4.7.   *Let* $\gamma_i$ *be defined as in Lemma 2.1. If*

$$B_n(t) = \sum_{j=1}^{\lfloor nt \rfloor} \Delta F_j^2 (-1)^j,$$

*then* $(F, B_n) \Rightarrow (F, \kappa B)$, *where* $\kappa^2 = 4\pi^{-1} + 2\pi^{-1}\sum_{i=0}^{\infty}(-1)^i\gamma_i^2 > 0$ *and* $B$ *is a standard Brownian motion independent of* $F$.

Proof.   Let

$$Y_n(t) = \sum_{j=1}^{\lfloor nt \rfloor}(-1)^j(\Delta F_j^2 - \sigma_j^2)$$

and

$$A_n(t) = \sum_{j=1}^{\lfloor nt \rfloor} (-1)^j \sigma_j^2,$$

so that $B_n = Y_n + A_n$. Note that $Y_n$ is of the form (3.35) with $h_j(x) = (-1)^j(x^2 - 1)$. As in the proof of Proposition 4.4, $\kappa$ is well-defined and strictly positive. Using the methods in the proof of Proposition 4.5, we have that $(F, Y_n) \Rightarrow (F, \kappa B)$. To complete the proof, observe that

$$\sup_{0 \le s \le t} |A_n(s)| \le \sigma_{\lfloor nt \rfloor}^2 + \sum_{j=1}^{\lfloor nt/2 \rfloor} |\sigma_{2j}^2 - \sigma_{2j-1}^2|,$$

and by (2.3),

$$|\sigma_{2j}^2 - \sigma_{2j-1}^2| \le \frac{\Delta t}{(2j-1)^{3/2}}.$$

Hence, $A_n \to 0$ uniformly on compacts. $\square$

COROLLARY 4.8. *Let $\kappa$ be as in Proposition 4.7. Define*

$$B_n(t) = \kappa^{-1} \sum_{j=1}^{\lfloor nt/2 \rfloor} (\Delta F_{2j}^2 - \Delta F_{2j-1}^2)$$

*and*

$$I_n(t) = \sum_{j=1}^{\lfloor nt/2 \rfloor} F(t_{2j-1})(F(t_{2j}) - F(t_{2j-2})).$$

*Let $B$ be a standard Brownian motion independent of $F$, and define*

$$I = \frac{1}{2}F^2 - \frac{\kappa}{2}B.$$

*Then $(F, B_n, I_n) \Rightarrow (F, B, I)$.*

PROOF. Note that

$$\left| \kappa B_n(t) - \sum_{j=1}^{\lfloor nt \rfloor} \Delta F_j^2 (-1)^j \right| \le \Delta F_{\lfloor nt \rfloor}^2,$$

so that by Proposition 4.7, $(F, B_n) \Rightarrow (F, B)$. Also note that by (1.5),

$$F(t)^2 = 2I_n(t) + \kappa B_n(t) + \varepsilon_n(t),$$

where $\varepsilon_n(t) = F(t)^2 - F(t_{2N})^2$ and $N = \lfloor nt/2 \rfloor$. The conclusion of the corollary is now immediate since $\varepsilon_n(t) \to 0$ uniformly on compacts. $\square$

**Acknowledgments.** Part of this material appeared in my doctoral dissertation and I would like to thank my advisors, Chris Burdzy and Zhen-Qing Chen, for their helpful advice. Regarding the generalization of my dissertation work, I gratefully acknowledge the helpful ideas of Sona Zaveri Swanson. I would also like to thank Bruce Erickson, Yaozhong Hu and David Nualart for helpful discussions and feedback, as well as Tom Kurtz for all of his support and guidance. Special thanks go to two anonymous referees for their many helpful comments and suggestions. This work was done while supported by VIGRE, for which I thank the NSF, University of Washington and University of Wisconsin–Madison.

Department of Mathematics
University of Wisconsin–Madison
480 Lincoln Dr.
Madison, Wisconsin 53706–1388
USA
E-mail: swanson@math.wisc.edu
URL: http://www.swansonsite.com/W